\documentclass{kbsjrnl}

\usepackage{amssymb}

\newcommand{\la}{\langle}
\newcommand{\ra}{\rangle}
\newcommand{\tl}{\triangleleft}
\newcommand{\Aut}{{\rm Aut}}
\newcommand{\PAut}{{\rm PAut}}
\newcommand{\mod}{{\rm mod\ }}

\newcommand{\cl}{{\rm cl}}
\newcommand{\Stab}{{\rm Stab}}
\newcommand{\Ker}{{\rm Ker}}

\newcommand{\operatorname}[1]{\mathop{\rm#1}\nolimits}

\newcommand{\F}{\mathbb{F}}
\newcommand{\Z}{\mathbb{Z}}
\newcommand{\Fix}[2]{\operatorname{fix}_{{#1}}({#2})}
\newcommand{\AGL}{\operatorname{AGL}}

\newcommand{\Fq}{\mathbb{F}_q}
\newcommand{\Fp}{\mathbb{Z}_p}
\newcommand{\poly}{\mathcal{P}}
\newcommand{\polyirr}{\poly_{\rm irr}}
  \renewcommand{\P}{\mathbb{P}}

\newcommand{\PSL}{\operatorname{PSL}}
\newcommand{\PGL}{\operatorname{PGL}}
\newcommand{\GL}{\operatorname{GL}}
\newcommand{\PGammaL}{\operatorname{P\Gamma L}}

\newcommand{\BSset}[1]{\mathcal{H}^{(#1)}}

  \newcounter{Case}
  \newenvironment{Case}[1][\unskip]{\refstepcounter{Case}\em
  \medskip \noindent Case \theCase\ #1.\ }{\unskip\upshape}
  \renewcommand{\theCase}{\arabic{Case}}

\begin{document}



\journame{}
\volnumber{}
\issuenumber{}
\issuemonth{}
\volyear{}

\received{}\revised{}

\authorrunninghead{}
\titlerunninghead{}




\begin{article}

\title{Transitive permutation groups of prime-squared degree}

\authors{Edward Dobson}
\email{dobson@math.msstate.edu}

\affil{Department of Mathematics and Statistics\\
        Mississippi State University\\
        Mississippi State, MS 39762, USA}

\authors{Dave Witte}
\email{dwitte@math.okstate.edu}
\affil{Department of Mathematics\\
        Oklahoma State University\\
        Stillwater, OK 74078, USA}

\editor{}

\abstract{We explicitly determine all of the transitive groups of degree
$p^2$, $p$ a prime, whose Sylow $p$-subgroup is not isomorphic to the
wreath product $\Z_p\wr\Z_p$.  Furthermore, we provide a general description
of the transitive groups of degree $p^2$ whose Sylow $p$-subgroup is
isomorphic to $\Z_p\wr\Z_p$, and explicitly determine most of them.  As
applications, we solve the Cayley Isomorphism problem for Cayley objects of
an abelian group of order $p^2$, explicitly determine the full automorphism
group of Cayley graphs of abelian groups of order $p^2$, and find all
nonnormal Cayley graphs of order~$p^2$.}

\keywords{permutation group, Cayley graph, $p$-group}

\section{Introduction}

In 1901, Burnside \cite{burn1} proved the following theorem.

\begin{theorem}[Burnside, \cite{burn1}]\label{norsylp}
Let $G$ be a transitive group of prime degree.  Then either $G$ is doubly
transitive or $G$ contains a normal Sylow $p$-subgroup.
\end{theorem}

\noindent If $G$ is a transitive group of prime degree and has a normal
Sylow $p$-subgroup, then it is not difficult to show that $G$ is
permutation isomorphic to a subgroup of $\AGL(1,p)$.  Similarly, it is also
straightforward to show that if $G$ is a transitive group of prime degree,
then $G$ has a normal Sylow $p$-subgroup if and only if $G$ is solvable.

A well-known consequence of the classification of the finite simple groups
is that all doubly transitive groups are known \cite[Theorem 5.3]{Pjc},
and hence all doubly transitive groups of prime degree are known.

Combining these results yields the following well-known classification of all
transitive groups of prime degree.

\begin{definition}
  We use the following standard notation.
  \begin{itemize}
  \item $S_p$ and $A_p$, respectively, denote a symmetric group and an
alternating group of degree $p$,
  \item $\AGL(d,p) = \Z_p^d \rtimes \GL(d,p)$ denotes the group of affine
transformations of the $d$-dimensional vector space $\F_p^d$ over $\F_p$,
  \item $M_{11}$ and $M_{23}$ denote Mathieu groups,
  \item $\PSL(d,q)$ and $\PGL(d,q)$, respectively, denote a projective special
linear group and a projective general linear group over the field $\F_q$ of
$q$~elements, and
  \item $\PGammaL(d,q)$ denotes the semidirect product of $\PGL(d,q)$
with the group of Galois automorphisms of $\F_q$.
  \end{itemize}
  \end{definition}

\begin{theorem}[{\cite[Cor.~4.2]{Feit}}] \label{degreep}
  Suppose $H$ is a subgroup of~$S_p$ that
contains~$\Z_p$. Let $S$ be a minimal normal subgroup
of~$H$, and let $N = N_{S_p}(S)$, so $S$ is simple and $S
\le H \le N$. Then $N/S$ is cyclic, and either:
  \begin{enumerate}
  \item $S = \Z_p$, $N = \AGL(1,p)$, and $N/S \cong \Z_{p-1}$; or
  \item $S = A_p$, $N = S_p$, and $N/S \cong \Z_2$; or
  \item $p = 11$ and $S = H = N = \PSL(2,11)$; or
  \item $p = 11$ and $S = H = N = M_{11}$; or
  \item $p = 23$ and $S = H = N = M_{23}$; or
  \item $p = (r^{d^{m+1}}-1)/(r^{d^m}-1)$ for some prime~$r$
and natural numbers $d$ and~$m$, and we have $S =
\PSL(d,r^{d^m})$, $N = \PGammaL(d,r^{d^m})$, and $N/S \cong
\Z_m$.
  \end{enumerate}
  \end{theorem}

In this paper, we will begin the classification of all transitive groups of
degree $p^2$. Our starting point is Theorem \ref{main2} below (proved at the end of \S\ref{BurnsideSection}), which provides an analogue of Burnside's Theorem \ref{norsylp}. This allows us to determine all of the transitive permutation groups of degree $p^2$ that do not have Sylow $p$-subgroup isomorphic to the wreath product $\Z_p\wr\Z_p$ (see Theorem \ref{main3}; the proof appears at the end of~\S\ref{ElemAbelSection}).   Furthermore, Proposition~\ref{wreath-prop} below describes how to construct every imprimitive permutation group of degree~$p^2$ whose Sylow $p$-subgroup is isomorphic to $\Z_p\wr\Z_p$. (This proposition is proved at the beginning of~\S\ref{WreathSection}.) Unfortunately, this proposition does not provide a complete classification of these permutation groups, because, in some cases, we do not have an explicit description of the possible choices for~$K$ and~$\phi$ in the conclusion of the proposition. However,
Theorem~\ref{degreep} describes the possible choices for $H$ and~$L$, and, in most cases, Subsection~\ref{codes-mod-n} describes the possible choices
for~$K$, and Subsection~\ref{crossed-homo} describes the possible choices
for~$\phi$. This leads to a complete classification for most primes~$p$;
specifically, the classification is complete for any prime~$p$, such that $p\notin \{11,23\}$ and $p \neq (q^d - 1)/(q-1)$, for every prime-power~$q$ and natural number~$d$. The problems that remain are described in a remark at the end of~\S\ref{WreathSection}.

\begin{definition}(cf.~Definition~\ref{recursion}) \label{Pp-1'}
  Let $P_{p-1}'$ denote the unique subgroup of $S_{p^2}$ (up to conjugacy) having order $p^p$ and containing a transitive subgroup isomorphic of $\Z_p\times\Z_p$ (and therefore not containing a transitive cyclic subgroup; see Lemma~\ref{wreath}).
  \end{definition}

\begin{theorem}\label{main2}
Let $G$ be a transitive permutation group of degree $p^2$, $p$ a prime,
with Sylow $p$-subgroup $P$.  Then either
  \begin{enumerate}
  \item $G$ is doubly transitive; or
  \item $P\tl G$; or
  \item  $P$ is equivalent to either $\Z_p \times \Z_p$, $P_{p-1}'$,
or $\Z_p\wr\Z_p$.
  \end{enumerate}
  \end{theorem}

\begin{theorem}\label{main3}
  Let $G$ be a transitive group of degree $p^2$ such that a Sylow $p$-subgroup $P$ of $G$ is not isomorphic to $\Z_p\wr\Z_p$.  Then, after replacing $G$ by a conjugate, one of the following is true.
  \begin{enumerate}
  \item \label{main3-doubly}
  $G$ is doubly transitive, and either
  \begin{itemize}
  \item $G = A_{p^2}$ or $S_{p^2}$; or
  \item $\PSL(d,q)\le G\le \PGammaL(d,q)$, where $(q^n-1)/(q-1) = p^2$; or
  \item $\Z_p \times \Z_p \le G\le\AGL(2,p)$.
  \end{itemize}

\item \label{main3-primitive}
  $G$ is simply primitive, has an elementary abelian Sylow $p$-subgroup and
either
  \begin{itemize}
  \item $\Z_p \times \Z_p \le G\le\AGL(2,p)$; or
  \item $G$ has a transitive, imprimitive subgroup~$H$ of index $2$, such that $H \le S_p \times S_p$ {\rm(}so $H$ is described in Lemma~\ref{SpxSp}{\rm)},
  \end{itemize}

\item \label{main3-Burnside}
  $G$ is imprimitive, $P\not\cong \Z_p \times \Z_p$, $P\not\cong P_{p-1}'$, and $P\tl G$, so $P \le G \le N_{S_{p^2}}(P)$ {\rm(}and $N_{S_{p^2}}(P)$ is
described in Lemma \ref{Np2} or~\ref{normelemabel}{\rm)};

\item \label{main3-P1}
  $G$ is imprimitive, $P = \Z_p \times \Z_p$ and $G \le S_p \times S_p$
{\rm(}so $G$ is described in  Lemma~\ref{SpxSp}{\rm)}; or

\item \label{main3-LFix}
  $G$ is imprimitive, $P = P_{p-1}'$, and $G = L P$, where
  $\Z_p \times \Z_p \le L \le S_p \times \AGL(1,p)$
  {\rm(}so $L$ is described in Lemma \ref{SpxSp}{\rm)}.
  \end{enumerate}
  \end{theorem}

\begin{definition} ({\cite[p.~168]{Cart-Eil}})
  Let $H$ be a group and let $A$ be an $H$-module. (That is,
$A$ is an abelian group on which $H$ acts by automorphisms.  Also note that abelian groups, when viewed as modules, are written additively.)
A function $\phi \colon H \to A$ is a \emph{crossed
homomorphism} if, for every $h_1,h_2 \in H$, we have
  $$ \phi( h_1 h_2 ) = h_2^{-1} \cdot \phi(h_1) + \phi(h_2)
.$$
  (This is equivalent to the assertion that the function $H
\to H \ltimes A$ defined by $h \mapsto \bigl( h, \phi(h)
\bigr)$ is a homomorphism.)
  \end{definition}

\begin{proposition} \label{wreath-prop}
  Let
  \begin{enumerate}
  \item $p$ be a prime;
  \item \label{wreath-prop-HL}
  $H$ and $L$ be transitive subgroups of~$S_p$, such
that $L$ is simple;
  \item \label{wreath-prop-K}
  $K$ be an $H$-invariant subgroup of the direct product $(N_{S_p}((L))^p$ containing $L^p$;
  \item \label{wreath-prop-phi}
  $\phi \colon H \to N_{S_p}(L)^p/K$ be a crossed
homomorphism; and
  \item $G_{H,L,K,\phi} = \{\, (h,v) \in H \ltimes
N_{S_p}(L)^p : \phi(h) = vK \,\} \le S_p \wr S_p$.
  \end{enumerate}
  Then $G_{H,L,K,\phi}$ is a transitive, imprimitive
subgroup of~$S_{p^2}$, such that a Sylow $p$-subgroup of~$G$ is isomorphic to
$\Z_p \wr \Z_p$.

Conversely, if $G$ is a transitive, imprimitive permutation group of
degree~$p^2$, such that a Sylow $p$-subgroup of~$G$ is isomorphic to $\Z_p
\wr \Z_p$, then $G$ is equivalent to $G_{H,L,K,\phi}$, for some $H$, $L$,
$K$, and~$\phi$ as above.
  \end{proposition}

To some extent, our proofs follow the outline that was used to determine all transitive groups of prime degree.

In \S\ref{KnownSection}, we recall known results that provide a
classification of certain types of transitive permutation groups of degree
$p^2$, namely, doubly transitive groups, groups with elementary abelian Sylow $p$-subgroup, and simply primitive groups. (Recall that a permutation group is \emph{simply primitive} if it is primitive, but not doubly transitive.)

In \S\ref{BurnsideSection}, we extend Theorem \ref{norsylp} to transitive
groups of degree $p^2$.  It follows by \cite[Theorem 3.4$'$]{W} that every
transitive group of prime power degree contains a transitive Sylow
$p$-subgroup; in particular, every transitive group of degree $p^2$ contains a transitive Sylow $p$-subgroup.  We first show that there are exactly $2p - 1$ transitive $p$-subgroups of $S_{p^2}$ up to permutation isomorphism and explicitly determine them (see Theorem~\ref{main1}). We also calculate the normalizer of each of these $p$-subgroups (see Lemmas~\ref{Np2} and~\ref{normelemabel}).  Next, we prove Theorem~\ref{main2}, which extends Burnside's Theorem~\ref{norsylp} to transitive groups of degree  $p^2$; that is, it determines which of these $2p - 1$ $p$-subgroups $P$ have the property that if $G\le S_{p^2}$ with Sylow $p$-subgroup $P$, then either $P\tl G$ or $G$ is doubly transitive.  Happily, only three of the $2p - 1$ transitive $p$-subgroups of $S_{p^2}$ fail to have this property.

We are left with the problem of finding every imprimitive or simply primitive subgroup of $S_{p^2}$ whose Sylow $p$-subgroup is one of the three transitive $p$-subgroups of $S_{p^2}$ for which the extension of Burnside's Theorem mentioned above does not hold.  Two of these $p$-subgroups are $\Z_p^2$ and the group $P_{p-1}'$ (see Definition \ref{Pp-1'} or~\ref{recursion}), which can, in a natural way, be regarded as the ``dual" of $\Z_p^2$. These two $p$-subgroups are considered in \S\ref{ElemAbelSection}, and the remaining $p$-subgroup, $\Z_p\wr\Z_p$, is considered in \S\ref{WreathSection}. However, as explained in the comments before Definition \ref{Pp-1'}, our results on $\Z_p\wr\Z_p$ are not complete.

In \S\ref{ApplicationSection}, we prove some straightforward applications of the above results that are of interest to combinatorialists.

We remark that some of the intermediate results (as well as some of the
applications) in this paper are known, and will give appropriate references
as needed.

\acknowledgments{This research was partially supported by grants from the
National Science Foundation (DMS-9623256 and DMS-9801136).}

\section{Some known results} \label{KnownSection}

\subsection{Doubly transitive groups}

The doubly transitive groups of degree $p^2$ can be determined much as in
the case of degree~$p$.  Burnside \cite[pg. 202]{burn3} proved the following result.

\begin{theorem}[Burnside, \cite{burn3}]
The socle of a finite doubly transitive group is either a regular
elementary abelian $p$-group, or a nonregular nonabelian simple group.
  \end{theorem}

If $G$ is doubly transitive of degree~$p^2$, then it is not difficult to
show that the socle of~$G$ is abelian if and only if $G \le \AGL(2,p)$.
(Note that an elementary abelian group of order $p^2$ is isomorphic to
$\Z_p^2$. Also, we remark that the doubly transitive subgroups of $\AGL(d,p)$ have been determined \cite{Hupp, Her}, cf.\ \cite[proof of Thm.~5.3]{Pjc}).

The doubly transitive groups with nonabelian socle are listed in \cite[Table on p.~8]{Pjc}. (This result relies on the classification of finite simple groups.) By inspection of this list, we see that the only such doubly transitive groups of degree $p^2$ are as follows.

\begin{theorem}\label{dt}
Let $G$ be a doubly transitive group of degree $p^2$ with nonabelian
socle.  Then either $G = A_{p^2}$, or $G = S_{p^2}$, or $\PSL(d,q)\le G\le
\PGammaL(d,q)$.
  \end{theorem}

\subsection{Imprimitive groups with elementary abelian Sylow {\mit p}-subgroup}

  In \cite[Proposition B]{J}, Jones determined the imprimitive
permutation groups of degree $p^2$ whose Sylow $p$-subgroup is elementary
abelian of order $p^2$.

\begin{theorem}[Jones, \cite{J}]\label{cross}
  Let $G$ be an imprimitive permutation group of degree $p^2$, where $p$ is
prime, such that a Sylow $p$-subgroup of $H$ is elementary abelian of order
$p^2$.  Then $G\le S_p\times S_p$.
  \end{theorem}

The following simple lemma expresses the conclusion of Theorem~\ref{cross}
more concretely.

\begin{lemma} \label{SpxSp}
  Let $G$ be a transitive subgroup of $S_{p^2}$. The following are equivalent:
  \begin{enumerate}
  \item \label{SpxSp-SpxSp}
  $G \le S_p\times S_p$.
  \item \label{SpxSp-N(H)xN(K)}
  There are transitive subgroups  $H$ and $K$ of $S_p$, such that $H \times K
\le G \le N_{S_p}(H) \times N_{S_p}(K)$.
  \item \label{SpxSp-f}
  There are transitive subgroups  $H$ and $K$ of $S_p$, and a homomorphism
$f\colon H\to N_{S_p}(K)/K$, such that
  $G = \{(\sigma,\tau)\in H\times N_{S_p}(K):f(\sigma) = \tau K\}$.
  \end{enumerate}
  \end{lemma}

\begin{proof}
  ($\ref{SpxSp-SpxSp} \Rightarrow \ref{SpxSp-N(H)xN(K)}$)
  Let $H = G \cap (S_p \times 1)$ and $K = G \cap (1 \times S_p)$, so $H,K \tl
G$. then
  $G \le N_{S_p \times S_p}(H \times K) = N_{S_p}(H) \times
N_{S_p}(K)$.

  ($\ref{SpxSp-SpxSp} \Rightarrow \ref{SpxSp-f}$)
  Let $H$ be the image of~$G$ under the projection to the first factor, and
let $K = G \cap (1 \times S_p)$. Then $G \le H \times N_{S_p}(K)$. By
definition of~$K$, we have
  $(G/K) \cap \bigl[ 1 \times \bigl( N_{S_p}(K)/K \bigr) \bigr] = 1$,
  so $G/K$ is the graph of a well-defined homomorphism $f \colon H \to
N_{S_p}(K)/K$. The desired conclusion follows.

  ($\ref{SpxSp-N(H)xN(K)} \Rightarrow \ref{SpxSp-SpxSp}$) and
($\ref{SpxSp-f} \Rightarrow \ref{SpxSp-SpxSp}$)  are obvious.
  \end{proof}

\subsection{Simply primitive groups}

  The simply primitive groups of degree $p^2$ are given by the following
theorem of Wielandt \cite[Theorems 8.5 and 16.2]{WO}. (Recall that any
subgroup~$H$ as in part~(\ref{primitive-SpxSp}) of this result is described in
Lemma~\ref{SpxSp}.)

\begin{theorem}[Wielandt, \cite{WO}]\label{primitive}
  Let $G$ be a simply primitive permutation group of degree $p^2$, where $p$
is prime.  Then the Sylow $p$-subgroups of $G$ are elementary abelian of
order $p^2$, and either
  \begin{enumerate}
  \item $G$ has a unique elementary abelian Sylow $p$-subgroup, or
  \item \label{primitive-SpxSp}
  $G$ has an imprimitive subgroup $H$ of index $2$ {\rm(}and, from
Theorem~\ref{cross}, we have $H \le S_p \times S_p${\rm)}.
  \end{enumerate}
  \end{theorem}

\section{The Extension of Burnside's Theorem} \label{BurnsideSection}

In view of the results in \S\ref{KnownSection}, this section is mainly
concerned with imprimitive groups $G$ of degree $p^2$ whose Sylow
$p$-subgroups $P$ are not elementary abelian.  We begin in a slightly more
general context.

Let $G$ be a transitive permutation group of degree $mp$ acting on $\Z_m\times \Z_p$ that admits a complete block system ${\cal B}$ of $m$ blocks of cardinality $p$.  We may suppose without loss of generality that $G$ acts on $\Z_m\times\Z_p$ such that ${\cal B} = \{\{i\}\times\Z_p:i\in\Z_p\}$.  If $g\in G$, then $g$ permutes the $m$ blocks of ${\cal B}$  and hence induces a permutation in $S_m$ denoted $g/{\cal B}$.  We define $G/{\cal B} = \{g/{\cal B}:g\in G\}$.  Let $\Fix{G}{\cal B} = \{g\in G:g(B) = B{\rm \ for\ every\ }B\in{\cal B}\}$.  Assume that  $\Fix{G}{\cal B}\not = 1$ so that a Sylow $p$-subgroup~$P_0$ of $\Fix{G}{\cal B}$ is nontrivial.  
Then $P_0$ is contained in $\la
z_i:i\in\Z_m\ra$, where each $z_i$ is a $p$-cycle that permutes the elements of $\{i\}\times\Z_p$.  For $h\in P_0$, we then have that $h = \prod_{i = 0}^{m-1}z_i^{a_i}$, $a_i\in\Z_p$.  Define $v\colon P_0\to Z_p^m$ by $v(h) = (a_0,a_1,\ldots,a_{m-1})$.

\begin{lemma}\label{code}
The set $\{v(h):h\in P_0\}$ is a linear code of length $m$ over $\F_p$.
\end{lemma}

\begin{proof}  As a linear code of length $m$ over $\F_p$ is simply an
$m$-dimensional vector space over $\F_p$, we need only show that
$\{v(h):h\in P_0\}$ is a vector space.   Note that for $g,h\in P_0$ and
$r\in\Z_p$, we have
  $$g^r = \left( \prod_{i=0}^{m-1}z_i^{a_i} \right)^r =
\prod_{i=0}^{m-1}z_i^{ra_i} $$
  and
  $$gh = \prod_{i=0}^{m-1}z_i^{a_i}\prod_{i=0}^{m-1}z_i^{b_i} =
\prod_{i=0}^{m-1}z_i^{a_i+b_i} .$$
  Hence $v(g^r) = rv(g)$ and  $v(gh) = v(g) + v(h)$, so $\{v(h):h\in P_0\}$
is a linear code. \end{proof}

\begin{definition} The code of Lemma \ref{code} will be denoted by $C_{\cal
B}$, and will be called the {\it code induced by ${\cal B}$}.  If $G$ admits a unique block system ${\cal B}$ of $m$ blocks of cardinality $p$, we say $C_{\cal B}$ is the code over $\F_p$ {\it induced} by $G$.  We remark that $C_{\cal B}$ depends upon the choice of the Sylow $p$-subgroup $P_0$, but that different choices of $P_0$ give monomially equivalent (that is, isomorphic) codes.
  \end{definition}

\begin{remark}
Lemma \ref{code} was proven in a less general context in \cite{Huff1}.
\end{remark}

\begin{lemma}\label{ccode}
If there exists $x\in G$ such that $x(i,j) = (i + 1,\alpha j + b_i)$,
$b_i\in\Z_p$, $\alpha\in\F_p^*$, then $\{v(h):h\in P_0\}$ is a cyclic
code of length $m$ over $\F_p$.  Conversely, if $C$ is a cyclic code
of length $m$ over $\F_p$, then there exists a group $G$ as above
such that $P_0 =
\bigl\{\, \prod_{i=0}^{m-1}z_i^{a_i}:(a_0,a_1,\ldots,a_{m-1})\in C \,\bigr\}$.
\end{lemma}

\begin{proof} From the form of~$x$, we know that $x$ normalizes
  $\la z_i : i \in \Z_m \ra$. also, because $x \in G$, we know that $x$
normalizes $\Fix{G}{\cal B}$. Thus, $x$ normalizes
  $\la z_i : i \in \Z_m \ra \cap \Fix{G}{\cal B} = P_0$.

For $h = \prod z_i^{a_i} \in P_0$ we have $x^{-1} h x = \prod z_i^{\alpha
a_{i+1}}$, so, because $x$ normalizes~$P_0$, we see that the linear code
$\{v(h):h\in P_0\}$ is cyclic.

  Conversely, define $x\colon \Z_m\times\Z_p\to\Z_m\times\Z_p$ by $x(i,j) =
(i + 1,j)$. Then it is also straightforward to check that $G =
\{x^ig:i\in\Z_m,g\in C\}$ will do.
  \end{proof}

In the following we consider the case $m = p$, where $\Z_p\times\Z_p$ is identified with $\Z_{p^2}$ via $(a,b)\mapsto a + bp$.

\begin{definition}\label{recursion}
  Let $a_{i,j} = {i \choose j}(-1)^{i-j}$.  A straightforward calculation
will show that $a_{i,j - 1} = a_{i+1,j}+a_{i,j}$.  For $1\le i\le p$, let
  $$\gamma_i = z_0^{a_{p-i,0}}z_1^{a_{p-i,1}}\ldots z_{p-1}^ {a_{p-i,p-1}}
.$$
  Define $\tau\colon \Z_{p^2}\to\Z_{p^2}$ by
  $$\tau(i) = i + 1 (\mod p^2)$$
and $\rho_1,\rho_2\colon \Z_p^2\to\Z_p^2$ by
  $$ \mbox{$\rho_1(i,j) = (i,j + 1)$ and
$\rho_2(i,j) = (i+1,j)$.}$$
Using the above identification of $\Z_p\times\Z_p$ with $\Z_{p^2}$, we have, for example,
  $$z_i(a + bp) =
  \cases{
  a + bp & \mbox{if $a\not = i$} \cr
  a + (b + 1)p\ (\mod p^2) & \mbox{if $a = i$,}}
  $$
$\rho_1 = \Pi_{i=1}^{p-1}z_i$ and $\rho_2(a + bp) = (a + 1)\ (\mod p) + bp$.  Hence $\tau = z_{p-1}\rho_2$.
   Let
  $$ \mbox{$P_i = \la\tau,\gamma_i\ra$ and $P'_i =
\la\rho_1,\rho_2,\gamma_i\ra$,} $$
  for $1\le i\le p$.  We remark that $P_p = P_p'\cong\Z_p\wr\Z_p$.  There
are thus $2p - 1$ distinct groups $P_i$, $P_i'$, $1\le i\le p$.
  \end{definition}

\begin{theorem}\label{main1}
  Let $G$ be a transitive group of degree $p^2$ with Sylow $p$-subgroup $P$.
Let $\vert P\vert = p^{i+1}$, $i\ge 1$.
  \begin{itemize}
  \item If $\tau\in P$, then $P = P_i$.
  \item If $\la\rho_1,\rho_2\ra\le P$, then $P = \alpha^{-1}P_i'\alpha$ for
some $\alpha\in\Aut(\Z_p^2)$.
  \end{itemize}
  \end{theorem}

\begin{proof} By \cite{MS}, if $C$ is a cyclic code of length $p$ over
$\F_p$, then $C$ has generator polynomial $f(x)$, where $f(x)$ divides $x^p
- 1$ in $\Z_p[x]$.  By the Freshman's Dream \cite{Hun}, $x^p -1 = (x - 1)^p$
so that $f(x) = (x - 1)^i$ for some $0\le i\le p - 1$.  As $C$ is generated by
the cyclic shifts of the vector $(a_{i,0},a_{i,1},\ldots,a_{i,p-1})$ where
$(x-1)^i = \sum_{j=0}^i a_{i,j}x^j$, we have $a_{i,j} = {i \choose
j}(-1)^{i-j}$.  Finally, we remark that the dimension of the code $C$ is $p -
i$.

Let $G$ be a transitive group of degree $p^2$ such that $\tau\in G$.  Let
$P$ be the Sylow $p$-subgroup of $G$ that contains $\tau$.  Then $P$ admits
a complete block system ${\cal B}$ of $p$ blocks of cardinality $p$ formed
by the orbits of $\la\tau^p\ra$.  Then $\vert\Fix{G}{\cal B}\vert = p^i$
and $C = \{ v(g):g\in\Fix{G}{\cal B}\}$ is a cyclic code of length $p$
over $\F_p$ by Lemma \ref{ccode}, so that $C$ contains $p^i$ codewords and
is thus of dimension $i$.  We conclude that $P =
\la\tau,\prod_{i=0}^{p-1}z_i^{a_i}:(a_0,a_1,\ldots,a_{p-1})\in C\ra =
\la\tau,\gamma_i\ra$.

If $\la\rho_1,\rho_2\ra\le P$, then again $P$ admits a complete
block system ${\cal B}$ formed by the orbits of $\la\delta\ra$, where
$\la\delta\ra = \la\rho_1\ra$ or $\la\rho_1^i\rho_2\ra$, $0\le i\le p
- 1$.  Hence there exists a group automorphism $\alpha$ of $\Z_p^2$
such that $\alpha^{-1}\delta\alpha = \rho_2$.  It then follows by
arguments above that $P =
\alpha^{-1}\la\rho_1,\rho_2,\gamma_i\ra\alpha$.
\end{proof}

\begin{remark}
  Every transitive group $G$ of degree $p^2$ contains a subgroup isomorphic
to either $\la\tau\ra$ or $\la\rho_1,\rho_2\ra$,
  as every Sylow $p$-subgroup of $G$ is transitive and contains a nontrivial
center.  Hence the above result determines all transitive $p$-subgroups of
$S_{p^2}$ up to isomorphism. \end{remark}

\begin{remark}
  Theorem \ref{main1} was already proven in \cite{HJP} for the case where $P$
contains a regular subgroup isomorphic to $\la\tau\ra$.
  \end{remark}

\begin{lemma}\label{wreath}
  Let $P$ be a transitive $p$-subgroup of $S_{p^2}$. Then $P$ admits a
complete block system ${\cal B}$ of $p$ blocks of cardinality $p$.

Furthermore, the following are equivalent:
  \begin{enumerate}
  \item \label{wreath-notboth}
  $P$ does not contain regular copies of both $\Z_{p^2}$ and $\Z_p^2$.
  \item \label{wreath-notwreath}
  $P\not\cong\Z_p\wr\Z_p$.
  \item \label{wreath-code}
  Letting $C$ be the code induced by ${\cal B}$, we have
$\sum_{i=0}^{p-1}a_i\equiv 0\ (\mod p)$, for every
$(a_0,a_1,\ldots,a_{p-1})\in C$.
  \end{enumerate}
  \end{lemma}

\begin{proof}
  ($\ref{wreath-notboth} \Rightarrow \ref{wreath-notwreath}$)
  If $P\cong\Z_p\wr\Z_p$, then $P$ is a Sylow $p$-subgroup of $S_{p^2}$, so
it is clear that $P$ contains both a regular subgroup isomorphic to
$\Z_{p^2}$ and a regular subgroup isomorphic to $\Z_p^2$.

  ($\ref{wreath-notwreath} \Rightarrow \ref{wreath-notboth}$)
  Assume $P$ contains regular copies of both $\Z_{p^2}$ and $\Z_p^2$.
Without loss of generality assume that $\tau\in P$.  As $P$ contains a
nontrivial center, $\tau^p\in Z(P)$, so that $P$ admits a complete block
system ${\cal B}$ of $p$ blocks of cardinality $p$ formed by the orbits of
$\la\tau^p\ra$.  As $P$ contains a regular subgroup isomorphic to $\Z_p^2$,
there exists $\tau_1,\tau_2\in P$ such that
$\la\tau_1,\tau_2\ra\cong\Z_p^2$.  As $\vert P/{\cal B}\vert = p$, we
assume without loss of generality that $\tau_2/{\cal B} = 1$ so that
$\vert\tau_1/{\cal B}\vert = p$.  As $\vert\tau_1\vert = p$,
$\tau_1^{-1}\tau^p\tau_1 = \tau^p$ and $\vert\tau^p\vert = p$, we may
assume that $\tau_2 = \tau^p$.  We regard $\Z_{p^2}$ as $\Z_p^2$.  Hence
$\tau(i,j) = (i+1,\delta_i(j))$ where $\delta_i(j) = j$, $0\le i\le p - 2$
and $\delta_{p-1}(j) = j + 1$.  Further, $\tau_1(i,j) = (i+ r,j + b_i)$,
$r,b_i\in\Z_p$.  As $\vert\tau_1\vert = p$, $\sum_{i=0}^{p-1}b_i\equiv 0\
(\mod p)$.  We assume without loss of generality that $r = 1$.  Then
$\tau^{-1}\tau_1(i,j) = (i,j + c_i)$ where $\sum_{i=0}^{p-1}c_i\equiv -1\
(\mod p)$.  Then $\Fix{P}{\cal B}  = \la\tau,\tau^{-j}\gamma_i\tau^j:1\le
j\le p - 1\ra$, for some $1\le i\le p$.  If $1\le k\le p - 1$ and $\psi\in
P_k$ with $\psi(i,j + d_i)$, we have that $\sum_{i=0}^{p-1}d_i\equiv 0\
(\mod p)$.  Hence $i = p$ and $\la\tau,\tau^{-1}\tau_1\ra\cong \Z_p\wr\Z_p$
as required.

  ($\ref{wreath-code} \Rightarrow \ref{wreath-notwreath}$)
  Obvious.

  ($\ref{wreath-notwreath} \Rightarrow \ref{wreath-code}$)
  If $P\not\cong\Z_p\wr\Z_p$, then, from the proof of Theorem~\ref{main1}, we
see that the generating polynomial of~$C$ is divisible by $x-1$. The
desired conclusion follows.
  \end{proof}

We now calculate the normalizers of $P_i $ and $P_i '$, $i\le p$. (We
remark that the normalizer of each $P_i$ was calculated in \cite{HJP}.)

\begin{definition}
  For $\beta\in\F_p^*$, define $\bar{\beta},\tilde{\beta}\colon
\Z_p^2\to\Z_p^2$ by
  $$ \mbox{ $\bar{\beta}(i,j) = (\beta i,j)$ and $\tilde{\beta}(i,j) = (i,\beta
j)$.}$$ For $\beta\in\Z_{p^2}^*$ define
  $\hat{\beta}\colon \Z_{p^2}\to\Z_{p^2}$ by
  $$\hat{\beta}(i) = \beta i .$$
  \end{definition}

\begin{remark}
Of course, $N_{S_{p^2}}(P_i)$ admits a (unique) complete block system
${\cal B}$ of $p$ blocks of cardinality $p$ formed by the orbits of
$\la\tau^p\ra$.  It is straightforward to show that
$\Fix{N_{S_{p^2}}(P_i)}{\cal B}$ is a $p$-group of order $p^i$, $1\le i\le p
- 1$.
  \end{remark}

\begin{lemma}[{\cite{HJP}}]\label{Np2}
  We have
  $$ N_{S_{p^2}}(P_i)
  =
  \cases{
  P_{i+1} \rtimes \bigl\{\, \hat\beta : \mbox{$\beta \in \Z_{p^2}^*$,
$\vert \beta \vert \in \{1, p-1\}$} \,\bigr\}
  & \mbox{if $1 \le i \le p-1$} \cr
  P_p \rtimes \{\,
\bar{\beta},\tilde{\beta}:\beta\in\F_p^*\,\}
   & \mbox{if $i = p$.}
  }$$
  \end{lemma}

\begin{proof}
  It is well known that
  $$ N_{S_{p^2}}(P_1) = N_{S_{p^2}}(\la\tau\ra)
  = \{\, x \mapsto ax + b:a\in\Z_{p^2}^*, b\in\Z_{p^2}\,\} $$ and it
essentially follows by arguments in \cite{APa} and was explicitly shown in
\cite{HJP} that
  $$ N_{S_{p^2}}(P_p) =  P_p \rtimes \{\,
\bar{\beta},\tilde{\beta}:\beta\in\F_p^*\,\} ,$$
  so we may assume $2 \le i \le p-1$.

We first show that $\vert N_{S_{p^2}}(P_i)\vert = (p-1)p^{i+2}$.

Let $X = \{\la x\ra:\la x\ra {\rm\ is\ a\ regular\ cyclic\ subgroup\ of\
}S_{p^2}\}$ and let $S_{p^2}$ act on $X$ by conjugation.  Denote the
resulting transitive permutation group on $X$ by $\Delta$.  Note that there
are $p^2!/(p-1)p^3 = [S_{p^2}:N_{S_{p^2}}(\la\tau\ra)]$ elements of $X$.
As $\la\tau^p\ra\le Z(P_i)$ and is the unique subgroup of $Z(P_i)$ of order
$p$, $N_{S_{p^2}}(P_i)$ admits a complete block system ${\cal B}$ of $p$
blocks of cardinality $p$, formed by the orbits of $\la\tau^p\ra$.  Observe
that if $\la x\ra\in X$ and $\la x\ra\le P_i$, then we may assume that $x =
\tau\gamma$, $\gamma\in\Fix{P}{\cal B}$.  Then $\vert\Fix{P}{\cal B}\vert
= p^i$, and by Lemma \ref{wreath} $\tau\gamma$ is a $p^2$-cycle for every
$\gamma\in\Fix{P}{\cal B}$ as every minimal transitive subgroup of $P_i$
is isomorphic to $\Z_{p^2}$.  Furthermore, there are exactly $p$ elements
of $\la\tau\gamma\ra$ contained in $\Fix{P}{\cal B}$.  We conclude that
$P_i$ contains $p^i/p = p^{i-1}$ elements of $X$.  Let $B =
\{\la\tau\gamma\ra:\gamma\in \Fix{P}{\cal B}\}$.  We first will show that
$B$ is a block of $\Delta$.

Let $\delta\in S_{p^2}$ such that $\delta^{-1}B\delta\cap B\not =
\emptyset$.  Then there exists $x = \tau\gamma$, $\gamma\in\Fix{P}{\cal B}$
such that $\delta^{-1}\la x\ra\delta\le P_i$.  Then $\delta^{-1}\la
x\ra\delta/{\cal B} = \la x\ra/{\cal B}$ and hence $\delta^{-1}\la
y\ra\delta/{\cal B} = \la\tau\ra/{\cal B}$ for every $\la y\ra\in B$.  Then
$\la\delta^{-1}B\delta\ra$ satisfies the hypothesis of Lemma \ref{ccode} (as
$\delta^{-1}\la x\ra\delta\le P_i$) so the code corresponding to
$\la\delta^{-1}B\delta\ra$ is the code corresponding to $\la B\ra$ which
implies $\la\delta^{-1}B\delta\ra = \la B\ra$ so that $\delta^{-1}B\delta =
B$ as required.  Hence the number of subgroups conjugate in $S_{p^2}$ to $P_i
= \la B\ra$ is the number of blocks conjugate to $B$ in $\Delta$.  As there
are $(p^2!/p^3(p-1))/p^{i-1}$ such blocks, $\vert N_{S_{p^2}}(P_i)\vert =
(p-1)p^{i+2}$ as required.

It is straightforward to check using the recursion formula given in
Definition \ref{recursion} that $\gamma_{i+1}\in N_{S_{p^2}}(P_i)$.  Note
that the result is clearly true for $i = 1$ as $\la\tau,\gamma_2\ra\le
N_{S_{p^2}}(\la\tau\ra)$.  Hence $N_{S_{p^2}}(\la\tau\ra)\le
N_{S_{p^2}}(\la\tau,\gamma_3\ra)$ as $\la\tau,\gamma_2\ra$ is the unique
Sylow $p$-subgroup of $N_{S_{p^2}}(\la\tau\ra)$.  Continuing inductively,
we have that $N_{S_{p^2}}(\la\tau\ra)\le N_{S_{p^2}}(\la\tau,\gamma_i)$ and
as $\vert\la N_{S_{p^2}}(\la\tau\ra),\gamma_i\ra\vert = (p-1)p^{i+2}$, the
result follows.
  \end{proof}

\begin{lemma}\label{normelemabel}
  We have
  $$ N_{S_{p^2}}(P_i') =
  \cases{
  \AGL(2,p) & \mbox{if $i = 1$} \cr
  P_{i+1}' \rtimes \{\, \bar{\beta},\tilde{\beta}: \beta\in\F_p^*\,\} &
\mbox{if $2\le i\le p - 1$} \cr
  P_p' \rtimes \{\, \bar{\beta},\tilde{\beta}: \beta\in\F_p^*\,\} &
\mbox{if $i = p$.}
  }$$
  \end{lemma}

\begin{proof}
  The case $i = 1$ is well known, and the case $i = p$ appears in Lemma
\ref{Np2}, so we may assume $2 \le i \le p - 1$.

Because $i \ge 2$, we know that $\la\tau_2\ra = Z(P_i')$ is characteristic in
$P_i'$, so that $\la\tau_2\ra\tl N_{S_{p^2}}(P_i')$.  Hence if $\delta\in
N_{S_{p^2}}(P_i')$, then $\delta(i,j) = (\sigma(i),\gamma j + b_i)$,
$\sigma\in S_p$, $\gamma\in\F_p^*$, $b_i\in\Z_p$.  Furthermore, $\sigma(i)
= \beta i + b$, $\beta\in\F_p^*$, $b\in\Z_p$.  It is
straightforward to check that $\bar{\beta}, \tilde{\gamma}\in
N_{S_{p^2}}(P_i')$, so we may assume $\beta = \gamma = 1$.  As $\tau_1\in
P_i'$, we may assume without loss of generality that $b = 0$.  As
$\gamma_1,\gamma_2,\ldots, \gamma_i\in P_i'$, we may assume, for $0\le k
\le i - 1$, that $b_k = 0$.  It is then straightforward to show that
$\delta\in\la\gamma_{i+1}\ra$.
  \end{proof}

\begin{definition}
  A code $C$ of length $m$ over $\F_p$ is said to be {\it degenerate} if
there exists $k\vert m$, $k\not = m$, and a code $D$ of length $k$ over
$\F_p$ such that $C = \oplus_{i=1}^{m/k}D$.  That is, a code $C$ is degenerate if $C$ consists $m/k$ repetitions of $D$.  If $C$ is not degenerate, we say that it is {\it nondegenerate}.
  \end{definition}

\begin{lemma}\label{killdt}
  Let $G\le S_{p^2}$ admit a complete block system ${\cal B}$ of $p$ blocks
of cardinality $p$.  If $\Fix{G}{\cal B}$ contains at least two Sylow
$p$-subgroups and $C_{\cal B}$ is nondegenerate, then  a Sylow $p$-subgroup
of $G$ is isomorphic to $\Z_p\times\Z_p$.
  \end{lemma}

\begin{proof}  We assume that $\la\tau\ra\le G$ or $\la\rho_1,\rho_2\ra\le
G$.  Let $P$ be a Sylow $p$-subgroup of $G$ that contains $\la\tau\ra$ or
$\la\rho_1,\rho_2\ra$.  Observe that $P$ cannot contain both $\la\tau\ra$
and $\la\rho_1,\rho_2\ra$, for then Lemma \ref{wreath} would imply
$P\cong\Z_p\wr\Z_p$, in which case $C_{{\cal B}}$ is degenerate.  If
$\Fix{G}{\cal B}$ contains at least two Sylow $p$-subgroups, then
$\Fix{G}{\cal B}\vert_B$ contains at least two Sylow $p$-subgroups for
every $B\in{\cal B}$ and hence, by the comments following Theorem
\ref{norsylp}, is nonsolvable.  By Theorem \ref{norsylp} $\Fix{G}{\cal
B}\vert_B$ is doubly transitive for every $B\in{\cal B}$.

Suppose, for the moment, that $\Fix{P}{\cal
B}$ is faithful on some block of ${\cal B}$. Then a Sylow $p$-subgroup of
$G$ has order $p^2$.  If $P = \la\rho_1,\rho_2\ra$, we are finished, so we
assume that $P = \la\tau\ra$ and hence $\Fix{P}{\cal B} = \la\tau^p\ra$.
Clearly $\la\tau^p\ra\vert_B$ is a Sylow $p$-subgroup of $\Fix{G}{\cal
B}\vert_B$, and if $N_{\Fix{G}{\cal B}\vert_B}(\la\tau^p\ra\vert_B) =
\la\tau^p\ra\vert_B$, then by Burnside's Transfer Theorem \cite[Theorem
4.3, pg. 252]{Gorenstein}, $\la\tau^p\ra\vert_B$ has a normal
$p$-complement in $\Fix{G}{\cal B}\vert_B$.  Whence $\Fix{G}{\cal
B}\vert_B$ admits a complete block system of $p$ blocks of cardinality $m$,
where $m\not = 1$, a contradiction.  Thus $N_{\Fix{G}{\cal
B}\vert_B}(\la\tau^p\ra\vert_B)\not = \la\tau^p\ra\vert_B$, so there exists
$\delta\in\Fix{G}{\cal B} - \la\tau^p\ra$ such that
$\delta^{-1}\tau^p\delta = \tau^{ap}$, $a\not = 1$.  By the remark
preceding Lemma \ref{Np2}, $\delta\not\in N_{S_{p^2}}(\la\tau\ra)$, so that
$\delta^{-1}\tau\delta\tau^{-1}\in\Fix{G}{\cal B}$, but
$\delta^{-1}\tau\delta\tau^{-1}\not\in\la\tau^p\ra$.  A straightforward
computation will then show that $\delta^{-1}\tau\delta\tau^{-1}$
centralizes $\la\tau^p\ra$.  As
$\delta^{-1}\tau\delta\tau^{-1}\in\Fix{G}{\cal B}$,
$\delta^{-1}\tau\delta\tau^{-1}\vert_B$ centralizes $\la\tau^p\ra\vert_B$,
and of course, $\la\tau^p\ra\vert_B$ is regular and abelian.  As a regular
abelian group is self-centralizing \cite[Proposition 4.4]{W}, we conclude
that $\delta^{-1}\tau\delta\tau^{-1}\vert_B\in\la\tau^p\ra\vert_B$.  Whence
$\delta^{-1}\tau\delta\tau^{-1}$ has order $p$ and
$\la\tau^p,\delta^{-1}\tau\delta\tau^{-1}\ra\le\Fix{G}{\cal B}$ and has
order $p^2$, a contradiction.

Henceforth, we assume that $\Fix{P}{\cal B}$ is not faithful on any block
of ${\cal B}$, so the Sylow $p$-subgroups of $\Fix{G}{\cal B}$ have order
at least $p^2$.
  Let $\gamma\in\Fix{G}{\cal B}$ such that $\gamma\vert_B\not = 1$ for the
fewest number of blocks $B\in{\cal B}$, and let ${\cal C} =
\{B\in{\cal B}:\gamma\vert_B\not = 1\}$.  If $\cup{\cal C}$ is a block
of $G$, then $C_{\cal B}$ is degenerate and we are finished.
Otherwise, define $\pi\colon \Fix{G}{\cal B}\to S_{\cup ({\cal B}-{\cal
C})}$ by $\pi(g) = g\vert_{\cup({\cal B}-{\cal C})}$.  Then
$\gamma\in\Ker(\pi)$.  As $\Ker(\pi)\tl\Fix{G}{\cal B}$ and
$\Fix{G}{\cal B}\vert_B$ is primitive for every $B\in{\cal B}$,
$\Ker(\pi)\vert_B$ is transitive for every $B\in{\cal C}$.  Hence we
may assume $\vert\gamma\vert = p$.  As any two Sylow $p$-subgroups of
$\Fix{G}{\cal B}$ are conjugate and one Sylow $p$-subgroup of
$\Fix{G}{\cal B}$ is contained in $\la z_i:i\in\Z_p\ra$, we assume without
loss of generality that $\gamma\in\la z_i:i\in\Z_p\ra$.  Finally,
observe that as $\gamma\in\Fix{G}{\cal B}$ such that $\gamma\vert_B\not =
1$ for the fewest number of blocks of ${\cal B}$, $\la\gamma\ra$ is a
Sylow $p$-subgroup of $\Ker(\pi)$.  As $\Ker(\pi)\tl\Fix{G}{\cal B}$,
$\Ker(\pi)\vert_B$ contains at least $2$ Sylow $p$-subgroups.  It then
follows by Burnside's Transfer Theorem that $N_{\Ker(\pi)}(\la\gamma\ra)\not =
\la\gamma\ra$.

Let $\gamma = z_{i_1}^{a_{i_1}}z_{i_2}^{a_{i_2}}\ldots
z_{i_r}^{a_{i_r}}$.  Let $\delta\in N_{\Ker(\pi)}(\la\gamma\ra)$ such
that $\delta\not\in\la\gamma\ra$.  Then $\delta^{-1}\gamma\delta =
z_{i_1}^{ba_{i_1}}z_{i_2}^{ba_{i_2}}\ldots z_{i_r}^{ba_{i_r}}$, for
some $b\in\Z_p^*$.  As $\cup{\cal C}$ is not a block of $G$, there
exists $\iota\in G$ such that $\iota^{-1}(\cup{\cal
C})\iota\cap(\cup{\cal C})\not = \emptyset$ and $\iota^{-1}(\cup{\cal
C})\iota\not = \cup{\cal C}$.  Let $\iota^{-1}\gamma\iota =
z_{j_1}^{c_{j_1}}z_{j_2}^{c_{j_2}}\cdots z_{j_r}^{c_{j_r}}$, for some
$j_1,j_2,\ldots,j_r\in\Z_{p^{k-1}}$ and $c_{j_\ell}\in\Z_p^*$.  Then
$\delta^{-1}\iota^{-1}\gamma\iota\delta(\iota^{-1}\gamma\iota)^{-b}\vert_B\not
= 1$ for fewer blocks of ${\cal B}$ than $\gamma$, a contradiction.
\end{proof}

\begin{definition}
  For a code $C$ of length $n$ over a field $\F_p$ of prime order $p$, let
$\Aut(C)$ be the group of all linear bijections of $K^n$ which map each
codeword of $C$ to a codeword of $C$ of the same weight.  Thus $\Aut(C)$ is
the subgroup of $M_n(\F_p)$ that map each codeword of $C$ to a codeword of
$C$, where $M_n(\F_p)$ is the set of all $n\times n$ monomial matrices
over $\F_p$.  That is, matrices with exactly one nonzero entry from
$\F_p$ in each row and column.
  \end{definition}

Let $[m_{ij}] = M\in M_n(\F_p)$.  Then $M = PD$, where $P = [p_{ij}]$ is the
permutation matrix given by $p_{ij} = 1$ if $m_{ij}\not = 0$ and $p_{ij} = 0$
otherwise and $D = [d_{ij}]$ is a diagonal matrix with $d_{ii} = m_{ij}$ if
$m_{ij}\not = 0$ and $d_{ij} = 0$ if $i\not = j$.  As the group of all
permutation matrices is simply the symmetric group on the coordinates of a
vector in $\F_p^n$, there is thus a canonical isomorphism between
$M_n(\F_p)$ and $S_n\ltimes (\F_p^*)^n$, with multiplication in $S_n\ltimes
(\F_p^*)^n$ given by $(\sigma,a)(\tau,b) = (\sigma\tau,(\sigma^{-1}b)a)$ and
$S_n\ltimes (\F_p^*)^n)$ acts on $\F_p^n$ by $(\sigma,d)(x) = \sigma(xd)$.
We will abuse notation and write that $(\sigma,d) = M\in M_n(F_p)$.  If
$(\sigma,d)\in\Aut(C)$, then $(\sigma,d)$ is {\it diagonal} if and only if
$\sigma = 1$.  Finally, we let $\PAut(C) = \{\sigma:(\sigma,d)\in\Aut(C)\}$.

\begin{theorem}[\cite{KS}, Theorem 1.3]\label{multipliers}
  If $C$ is a nontrivial code such that $\PAut(C)$ is primitive, then $C$ is
nondegenerate and every diagonal automorphism of $C$ is scalar.
  \end{theorem}

\begin{definition}
  A code $C$ of length $p$ over $\F_p$ is {\it affine invariant} if
$\AGL(1,p)\le \PAut(C)$.
  \end{definition}

Let $\Pi\le S_{p^2}$ be a transitive $p$-group.  Then $\Pi$ admits a complete
block system ${\cal B}$ of $p$ blocks of cardinality $p$, formed by the
orbits of a semiregular element of order $p$ contained in the center of
$\Pi$.  By Theorem \ref{main1}, we may assume that $\Pi = \Pi_i$ or
$\Pi_i'$.  Note that if $G\le S_{p^2}$ admits ${\cal B}$ as a complete block
system with Sylow $p$-subgroup $\Pi$, then conjugation by an element of $G$
induces an automorphism of $C_{\cal B}$.  It then follows by Lemmas \ref{Np2}
and \ref{normelemabel} that $C_{\cal B}$ is affine invariant.

\begin{theorem}[\cite{Huff}, Theorem 3.17]\label{bigtool}
  Let $C$ be a cyclic code of length $p$ over $\Z_p$ such that $\vert
C\vert\not = p$, $p^{p-1}$ or $p^p$.  Then $\Aut(C) =
\{(\sigma,d):\sigma\in\AGL(1,p){\rm\ and\ }d\in\F_p^*\}$.
  \end{theorem}

\begin{lemma}\label{topnotdt}
Let $G\le S_{p^2}$ be transitive such that $G$ admits a complete block
system ${\cal B}$ of $p$ blocks of cardinality $p$.  Furthermore,
assume either $\Pi_1$ or $\Pi_{p-2}$ is a Sylow $p$-subgroup of $G$ and $\Fix{G}{\cal B}$ is a $p$-group.  Then $G/{\cal B}\le\AGL(1,p)$.
\end{lemma}

\begin{proof} We argue by contradiction and thus assume that $G/{\cal B}$ has at least $2$ Sylow $p$-subgroups and is doubly transitive.  As $G/{\cal B}$ has at least $2$ Sylow $p$-subgroups and a Sylow $p$-subgroup of $G$ is $\Pi_1$ or $\Pi_{p-2}$, there exists $g\in G$
such that $\vert g/{\cal B}\vert = p$, $g$ is of order a power of
$p$, but $g/{\cal B}\not\in\la\tau\ra/{\cal B}$.  As $g$ is contained in a
Sylow $p$-subgroup of $G$ (which is isomorphic to $\Pi_1$ or $\Pi_{p-2}$), we have that $\la g\ra$ is cyclic and regular of order $p^2$.  As $G$ admits ${\cal B}$ as a complete block system and $\Fix{G}{\cal B}$ is a $p$-group, we have that $g(i + jp) = \sigma(i) + (\alpha_ij +
b_i)p$, $\sigma\in S_p$, $\alpha_i\in\F_p^*$, $b_i\in\Z_p$.  It follows by 
\cite[Theorem 3.16 (i)]{Huff} that $\alpha_i = \alpha_j$ for all $i,j\in\Z_p$.  As $\vert g\vert = p^2$, we have that $\alpha_i = 1$ for all $i,j\in\Z_p$.  As $\vert g\vert = p^2$, $\sum_{i=0}^{p-1}b_i\equiv c\not\equiv 0\ (\mod p)$.  Then there exists $r\in\Z_p$ such that $g\tau^r/{\cal B}$ has a fixed point, and $g\tau^r/{\cal B}\not = 1$.  Let $d = \vert g\tau^r/{\cal B}\vert$, and ${\cal O}_1,{\cal O}_2,\ldots,{\cal O}_s$ the orbits of $\la g\tau^r/{\cal B}\ra$.  Let $d_i = \vert{\cal O}_i\vert$.  Then $(g\tau^r)^d\in\Fix{G}{\cal B}$ so that if $(g\tau^r)^d =
z_0^{b_0}z_1^{b_1}\cdots z_{p-1}^{b_{p-1}}$, then
$\sum_{i=0}^{p-1}b_i\equiv 0\ (\mod p)$.  Let $g\tau^r(i + jp^{k-1})
= \delta(i) + (j + c_i)p$.  Then $\sum_{i=0}^{p-1}c_i\equiv c +
r\ (\mod p)$.  Further,

\begin{eqnarray*}
\sum_{i=0}^{p-1}b_i & = & \sum_{i=0}^s \left( \sum_{j\in{\cal
O}_i}a_i \right) \\
       & = & \sum_{i=0}^sd_i \left( \frac{d}{d_i}\sum_{j\in{\cal O}_i}c_i
\right)\\
       & = & \sum_{i=0}^{p-1}dc_i\\
       & \equiv & d(c+r)\ (\mod p).
\end{eqnarray*}

\noindent Note that as a Sylow $p$-subgroup of $G/{\cal B}$ is cyclic,
$\gcd(d,p) = 1$.  As $\sum_{i=0}^{p-1}b_i\equiv 0\ (\mod p)$ and
$(d,p) = 1$, we have that $c\equiv - r\ (\mod p)$.  However, there exists
$t\in\Z_p^*$, $t\not = r$ such that $g\tau^t/{\cal B}$ has a fixed
point, and analogous arguments will show that $c\equiv - t\ (\mod p)$, a
contradiction.
\end{proof}

Let $G$ be a transitive group acting on $\Omega$, and let $G$ act
on $\Omega\times\Omega$ by $g(\alpha,\beta) = (g(\alpha),g(\beta))$.  Let
${\cal O}_0,{\cal O}_1,\ldots,{\cal O}_n$ be the orbits of $G$ acting on
$\Omega\times\Omega$.  We call the orbit
$\{(\alpha,\alpha):\alpha\in\Omega\}$  the {\it trivial orbit}.  Assume
${\cal O}_0$ is the trivial orbit, and define directed graphs
$\Gamma_1,\ldots,\Gamma_n$ by $V(\Gamma_i) = \Omega$ and $E(\Gamma_i) = {\cal
O}_i$ for each $1\le i\le n$.  The graphs $\Gamma_i$ are {\it orbital
digraphs} of $G$.  Note that $G\le\Aut(\Gamma_i)$ for all $1\le i\le n$.  We
define the {\it 2-closure of $G$}, denoted $\cl(G)$ to be
$\cap_{i=1}^{n}\Aut(\Gamma_i)$.

\begin{proof}[{Proof of Theorem~\ref{main2}}]  Let $G\le S_{p^2}$ be
transitive such that $P\not\cong P_1'$, $P_{p-1}'$, or $\Z_p\wr\Z_p$.  If
$\tau\in P$, then, as $\Z_{p^2}$ is a Burnside group \cite{W}, $G$ is
either doubly transitive or imprimitive (this may also be obtained from a result of Guralnick \cite{Gur}), so we may assume $G$ is imprimitive.
If $\tau\not\in G$, then, as $p^3\vert\ \vert G\vert$, we see from
\cite{APa, G} that any orbital digraph $\Gamma_i$ of $G$ is either
complete, empty, or a nontrivial wreath product.  We conclude that
$\tau\in\cl(G)$ so that $G$ is again imprimitive or doubly transitive.  We
consider the cases $P\cong P_i$ and $P\cong P_i'$ for some $i$ separately.

If $P\cong P_i$, then by Lemma \ref{killdt} $\Fix{G}{\cal B}$ contains a
unique Sylow $p$-subgroup.  If $\Fix{G}{\cal B}$ is a $p$-group and $P\cong P_1$ or $P_{p-2}$, then it follows by Lemma \ref{topnotdt} that $G/{\cal B}\le AGL(1,p)$.  Then $G/{\cal B}$ contains exactly one Sylow $p$-subgroup, and $P\tl G$.  If $\Fix{G}{\cal B}$ is a $p$-group and $P\not\cong P_1$ or $P_{p-2}$, then it follows by Theorem \ref{bigtool} that $G/{\cal B}\le AGL(1,p)$ and the result follows.  If $\Fix{G}{\cal B}$ is not a $p$-group, then observe that conjugation of the unique Sylow $p$-subgroup of $\Fix{G}{\cal B}$ by an element of $G$ induces an automorphism of $C_{\cal B}$, the code induced by ${\cal B}$.  Let $g\in\Fix{G}{\cal B}$ such that $g$ is not in the unique Sylow $p$-subgroup of $G$.  Let $(\sigma,d)$ be the automorphism of $C_{\cal B}$ induced by conjugating the unique Sylow $p$-subgroup of $\Fix{G}{\cal B}$ by $g$.  Then $\sigma = 1$ so that $(\sigma,d)$ is diagonal.  As $C_{\cal B}$ has prime length, $\PAut(C_{\cal B})$ is primitive, and hence by Theorem \ref{multipliers}, we have that $(\sigma,d)$ is scalar.  Whence $g(i + jp) = i + (dj + b_i)p$.  A straightforward computation will then show that $g^{-1}\tau g\tau^{-1}$ is in the unique Sylow $p$-subgroup of $\Fix{G}{\cal B}$ so that $g^{-1}\tau g\in P_i$.  Thus $g\in N_{S_{p^2}}(P_i)$.  However, by Lemma \ref{Np2}, $\Fix{N_{S_{p^2}}(P_i)}{\cal B}$ contains no element of order relatively prime to $p$, a contradiction.

If $P\cong P_i'$, then by Lemma \ref{killdt}, $\Fix{G}{\cal B}$ contains
exactly one Sylow $p$-subgroup, and by arguments in the immediately
preceding paragraph, if $g\in\Fix{G}{\cal B}$ such that $\vert g\vert
\not = p$, then $g(i,j) = (i, d^{-1}j + b_i)$, $d\in\F_p^*$, $b_i\in\Z_p$.
It is then straightforward to verify that $g\in N_{S_{p^2}}(P_i')$.  By
Theorem \ref{bigtool}, $\la\tau_1\ra/{\cal B}\tl G/{\cal B}$.  Thus $H =
\la P_i',\Fix{G}{\cal B}\ra\tl G$, and $P_i'\tl H$.  As $P_i'$ is a
Sylow $p$-subgroup of $G$, it is the unique Sylow $p$-subgroup of $H$ and
so $P_i'$ is characteristic in $H$.  As $H\tl G$, we have that $P_i'\tl
G$.
  \end{proof}

\section{Overgroups of $\Z_p^2$ and its ``dual"} \label{ElemAbelSection}

\begin{lemma}\label{elemabel-dual}
  Let $G$ be an imprimitive group of degree $p^2$ with Sylow $p$-subgroup
isomorphic to $P_{p-1}'$.  Let ${\cal B}$ be the unique complete block
system of $P_{p-1}'$ of $p$ blocks of cardinality $p$.  Then there exists a
subgroup $H\le S_p\times\AGL(1,p)$ such that $G$ is equivalent to
$H\cdot\Fix{P_{p-1}'}{\cal B}$.
  \end{lemma}

\begin{proof}  We assume without loss of generality that a Sylow
$p$-subgroup of $G$ is $P_{p-1}'$ so that $G$ admits ${\cal B}$ as a
complete block system.  As $P_{p-1}'$ is a Sylow $p$-subgroup of $G$,
$C_{\cal B}$ is nondegenerate.  It then follows by Lemma \ref{killdt} that
$\Fix{G}{\cal B}$ is solvable and hence $\Fix{G}{\cal
B}\vert_B\le\AGL(1,p)$ for every $B\in{\cal B}$.  We now show that if $g\in
G$, then $g\in N_{S_{p^2}}(\la\tau_2\ra)$.

  As $\Fix{G}{\cal B}\vert_B\le\AGL(1,p)$ for every $B\in{\cal B}$ if $g\in
G$, then $g(i,j) = (\sigma(i),\alpha_ij + b_i)$, where $\sigma\in S_p$,
$\alpha_i\in\F_p^*$ and $b_i\in\Z_p$.  To show that $g\in
N_{S_{p^2}}(\la\tau_2\ra)$ is suffices to show that $\alpha_i = \alpha_j$
for every $i,j\in\Z_p$.  As a Sylow $p$-subgroup of $G$ is $P_{p-1}'$, by
Lemma \ref{wreath} $C_{\cal B}$ consists of all vectors in $\Z_p^p$ such
that the sum of the coordinates is congruent to $0$ modulo $p$.  We
conclude that $z_iz_j^{p-1}\in\Fix{G}{\cal B}$ for every $i\not =
j\in\Z_p$.  Then $gz_iz_j^{p-1}g^{-1} =
z_{\sigma(i)}^{\alpha_i}z_{\sigma(j)}^{\alpha_j(p-1)}$.  As
$gz_iz_j^{p-1}g^{-1}\in\Fix{G}{\cal B}$, we must have that $\alpha_i +
\alpha_j(p-1)\equiv 0\ (\mod\ p)$, so that $\alpha_i = \alpha_j$.  Thus if
$g\in G$, then $g\in N_{S_{p^2}}(\la\tau_2\ra)$.

  Let $B\in{\cal B}$ and $L = \{g\in G:g(B) = B\}$.  Then $L$ has a unique
Sylow $p$-subgroup $P$, namely the Sylow $p$-subgroup of $\Fix{G}{\cal B}$,
so that $P\tl L$.  By the Schur-Zassenhaus Theorem \cite[Theorem 2.1, pg.
221]{Gorenstein} $L$ contains a $p'$-subgroup $M$ which is a complement to
$P$ in $L$.  Let $P' = \la z_i:i\in\Z_p\ra$ and $W = G\cdot P'$.  Note that
$\vert W\vert = p\cdot\vert G\vert$, and let $K = \{w\in W:w(B) = B\}$.
Then $P'\tl K$ and is the Sylow $p$-subgroup of $K$.  Again by the
Schur-Zassenhaus Theorem, as $P'$ is solvable, any two $p'$-subgroups of
$K$ are conjugate in $W$.  As $M$ is a $p'$-subgroup of $L$, $M$ is a
$p'$-subgroup of $K$.  Recall that if $g\in G$, then $g(i,j) =
(\sigma(i),\alpha j + b_i)$, $\sigma\in S_p$, $\alpha\in\F_p^*$, and
$b_i\in\Z_p$.  For $g\in G$, define $\hat{g}\colon S_{\Z_p^2}\to
S_{\Z_p^2}$ by $\hat{g}(i,j) = (\sigma(i),\alpha j)$, and let $\hat{M} =
\{\hat{g}:g\in M\}$.  Clearly $\hat{M}$ is a subgroup and $\vert
\hat{M}\vert = \vert M\vert$.  Furthermore, as $P' = \la z_i:i\in\Z_p\ra$,
$\hat{M}\le K$.  Thus $\hat{M}$ is also a $p'$ subgroup of $K$ so that
there exists $\delta\in P'$ such that $\delta^{-1}M\delta = \hat{M}$.  Let
$G' = \delta^{-1}G\delta$.  Clearly $\hat{M}\le G'$ and as $\delta\in P'\le
N_{S_{p^2}}(P_{p-1}')$ we have that $P_{p-1}'\le G'$.  Let $H =
\la\tau_1,\hat{M}\ra$.  Then $H\le G'$ and for every $h\in H$, $h(i,j) =
(\sigma(i),\alpha j)$ so that $H\le S_p \times\AGL(1,p)$.  As $\vert M\vert
= \vert G\vert/p^p$ we have that $\vert H\vert = \vert G\vert/p^{p-1}$ so
that $\vert H\cdot\Fix{P_{p-1}'}{\cal B}\vert = \vert G\vert = \vert
G'\vert$.  As $H\cdot\Fix{P_{p-1}'}{\cal B}\le G'$, we conclude that $G'
= H\cdot\Fix{P_{p-1}'}{\cal B}$ and the result follows.
  \end{proof}

\begin{proof}[{Proof of Theorem~\ref{main3}}]
   (\ref{main3-doubly}) follows from Lemma \ref{dt}, and
(\ref{main3-primitive}) follows from Theorem \ref{primitive}.

Thus, we assume, henceforth, that $G$ is imprimitive. By Theorem \ref{main1}
the Sylow $p$-subgroups of $G$ are isomorphic to $P_i$ or $P_i'$, $1\le i\le
p$.  If no Sylow $p$-subgroup of $G$ is isomorphic to $P_1'$ or $P_{p-1}'$,
then (\ref{main3-Burnside}) follows from Theorem \ref{main2}.  If a Sylow
$p$-subgroup of $G$ is isomorphic to $P_1'$, then (\ref{main3-P1}) follows
from Theorem \ref{cross} and Lemma \ref{SpxSp}.  Finally, if a
Sylow $p$-subgroup of $G$ is isomorphic to $P_{p-1}'$, then (\ref{main3-LFix})
follows from Lemma \ref{elemabel-dual}.
  \end{proof}

\section{Imprimitive subgroups that contain a Sylow
$p$-subgroup of~$S_{p^2}$} \label{WreathSection}

Note that $\Z_p \wr \Z_p$ is a Sylow $p$-subgroup of $S_{p^2}$.

\begin{proof}[{Proof of Proposition~\ref{wreath-prop}}]
  ($\Rightarrow$) Because $N_{S_p}(L)/L$ is cyclic
(see~Theorem \ref{degreep}), we know that $\bigl( N_{S_p}(L) /L
\bigr)^p$ is abelian, so it is obvious that $K/L^p$ is a
normal subgroup of $\bigl( N_{S_p}(L) /L \bigr)^p$; hence $K$
is a normal subgroup of $N_{S_p}(L)^p$.
  Then, because $\phi$ is a crossed homomorphism and $K$ is
$H$-invariant, it is easy to verify that $G_{H,L,K,\phi}$ is
closed under multiplication. Therefore, it is a subgroup of
$S_p \wr S_p$.

It is straightforward to verify that $K$ is a normal
subgroup of $G_{H,L,K,\phi}$, and we have  $G_{H,L,K,\phi}/K
\cong H$, so $|G_{H,L,K,\phi}|$ is divisible by $|K| |H|$.
Because $K \supset L^p$, this implies that
$|G_{H,L,K,\phi}|$ is divisible by $p^{p+1}$. Therefore,
$G_{H,L,K,\phi}$ contains a Sylow $p$-subgroup of $S_p \wr
S_p$, so $G_{H,L,K,\phi}$ is transitive. Because
$G_{H,L,K,\phi} \le S_p \wr S_p$, we know that
$G_{H,L,K,\phi}$ is imprimitive.

\medskip

  ($\Leftarrow$) Because $G$ is imprimitive, we may assume
that $G \le S_p \wr S_p$. Then, because $p^{p+1} \mid
|G|$, we know that $G$ contains a Sylow $p$-subgroup of $S_p
\wr S_p$; assume, without loss of generality, that $G$
contains $\Z_p \wr \Z_p$. In particular, $G$ admits a unique
block system~$\mathcal{B}$, consisting of $p$~blocks of
cardinality~$p$.

  Let $H = G/\mathcal{B} \le S_p$, let $K =
\Fix{G}{\cal B}$, and let $\hat L$ be the smallest
normal subgroup of~$G$ that contains $1 \wr \Z_p$. It is easy
to see that $\hat L = 1 \wr L \cong L^p$, for some
transitive, simple subgroup~$L$ of~$S_p$.

The map $g \mapsto g/\mathcal{B}$ is a homomorphism from~$G$
onto~$H$, with kernel~$K$. Thus, there is an isomorphism
$\hat\phi \colon H \to G/K$, given by $h =
\hat\phi(h)/\mathcal{B}$.
  Because $1 \wr L$ is normal in~$G$, we know that $G \le
H \wr N_{S_p}(L)$, so we may write $\hat\phi(h) = h \,
\phi(h)$, with $\phi(h) \in \bigl( 1 \wr N_{S_p}(L)^p
\bigr) /K$. Because $\hat\phi$ is a homomorphism, it is
straightforward to verify that $\phi$ is a crossed
homomorphism.
  \end{proof}

The assumption that $L$ is simple is not necessary in the
definition of $G_{H,L,K,\phi}$, but this restriction makes
$L$ unique (up to conjugacy). For a given group~$G$, the
corresponding $H,L,K,\phi$ are not uniquely determined, but
the following simple lemma describes how to tell whether
$G_{H_1,L_1,K_1,\phi_1}$ is equivalent to
$G_{H_2,L_2,K_2,\phi_2}$.

\begin{definition}({cf.~\cite[Prop.~4.1]{Cart-Eil}})
  Let $H$ be a group, let $A$ be an $H$-module, and let
$\phi_1,\phi_2 \colon H \to A$ be crossed homomorphism. We
say that $\phi_1$ is \emph{cohomologous} to~$\phi_2$ if
there is an element~$a$ of~$A$, such that, for every $h
\in H$, we have
  $$ \phi_1( h ) - \phi_2(h) = h^{-1} a - a .$$
  (This is equivalent to the assertion that the homomorphisms
$h \mapsto \bigl( h, \phi_1(h) \bigr)$ and $h \mapsto \bigl(
h, \phi_1(h) \bigr)$ are conjugate via an element of~$A$.)

We remark that the equivalence classes of this equivalence
relation are, by definition, the elements of the cohomology
group $H^1(H,A)$.
  \end{definition}

\begin{lemma}
  Let $H_i,L_i,K_i,\phi_i$ be as in
Proposition~\ref{wreath-prop}, for $i = 1,2$.
  \begin{enumerate}
  \item If $G_{H_1,L_1,K_1,\phi_1}$ is equivalent to
$G_{H_2,L_2,K_2,\phi_2}$, then $L_1$ is conjugate to~$L_2$
{\upshape(}in~$S_p${\upshape)}.
  \item  If $G_{H_1,L_1,K_1,\phi_1}$ is equivalent to
$G_{H_2,L_2,K_2,\phi_2}$, and $L_1 = L_2$,
then there exists $g \in S_p$, such that, letting $\hat g =
(g,1) \in S_p \wr S_p$, we have
  \begin{enumerate}
  \item $g H_1 g^{-1} = H_2$;
  \item $\hat g K_1 \hat g^{-1} = K_2$; and
  \item $\phi_1^g$ is cohomologous to~$\phi_2$,
  where $\phi_1^g \colon H_2 \to \bigl( N_{S_p}(L_2)/L_2
\bigr)^p $ is defined by
  $\phi_1^g(h) = \hat g \, \phi_1( g^{-1} h g ) \, \hat g^{-1}
$.
  \end{enumerate}
  \end{enumerate}
  \end{lemma}

\begin{proof}
  Let $h \in S_{p^2}$, with $h
G_{H_1,L_1,K_1,\phi_1} h^{-1} =
G_{H_2,L_2,K_2,\phi_2}$, and let $\mathcal{B}$ be the unique complete
block system for $\Z_p \wr \Z_p$. Because  $\Z_p \wr \Z_p$
is contained in both $G_{H_1,L_1,K_1,\phi_1}$ and
$G_{H_2,L_2,K_2,\phi_2}$, the uniqueness of~$\mathcal{B}$ implies
that $h \mathcal{B} = \mathcal{B}$; thus, $h \in S_p \wr S_p$, so we
may write
   $h = ( g, x )$, with $g \in S_p$ and $x \in (S_p)^p$.
  Because $L_1 = L_2$, we must have $x \in  N_{S_p}(L_2)^p$.

Because $N_{S_p}(L_2)/L_2$ is abelian, this
implies that $x$ normalizes~$K_2$, so we must have $g K_1
g^{-1} = K_2$.

Because $H_i = G_{H_i,L_i,K_i,\phi_i}/\mathcal{B}$, we must have $g
H_1 g^{-1} = H_2$.

Replacing $G_{H_1,L_1,K_1,\phi_1}$ by its conjugate
under $\hat g$, we may assume that $g = 1$, so $H_1
= H_2$, $K_1 = K_2$, and $\phi_1^g = \phi_1$. Because
  \begin{eqnarray*}
  x \, \bigl\{\, \bigl(h,\phi_1(h) \bigr) : h \in H_1
\,\bigr\} \, x^{-1}
  &=& \frac{x \, G_{H_1,L_1,K_1,\phi_1} \, x^{-1}}{L_1}
  = \frac{G_{H_2,L_2,K_2,\phi_2}}{L_2} \\
  &=&  \bigl\{\, \bigl(h,\phi_2(h) \bigr) : h \in H_2
\,\bigr\},
  \end{eqnarray*}
  we see that $\phi_1$ is cohomologous to~$\phi_2$.
  \end{proof}

\subsection{Cyclic codes modulo~{\mit n}}
  \label{codes-mod-n}
  The Chinese Remainder Theorem (Lemma~\ref{CRT-modn}) reduces
the study of codes modulo~$n$ to the case where $n$~is a
prime power. Assuming that $p \nmid n$, the problem can often
be further reduced to the case where $n$ is prime
(see Lemma~\ref{reducetoprime} and
Remark~\ref{goodp-K-notPSL}). This reduced case is
considered in Lemma~\ref{modp-allgroups}.

\begin{lemma}
[{cf.~\cite[Thm.~1.2.13, p.~8]{Gorenstein}}]\label{CRT-modn}
  Let $n = n_1 n_2 \cdots n_r$, where each $n_i$ is a prime
power, and $n_1, n_2, \cdots, n_r$ are pairwise relatively
prime.
  \begin{enumerate}
  \item We have
  $(\Z_n)^p \cong (\Z_{n_1})^p \oplus
(\Z_{n_2})^p \oplus \cdots \oplus (\Z_{n_r})^p$.
  \item For any subgroup~$C$ of~$(\Z_{n_1})^p \oplus
(\Z_{n_2})^p \oplus \cdots \oplus (\Z_{n_r})^p$,
we have
  $$ C =
  \bigl( C \cap (\Z_{n_1})^p \bigr) \oplus
  \bigl( C \cap (\Z_{n_2})^p \bigr) \oplus
  \cdots
  \bigl( C \cap (\Z_{n_r})^p \bigr)
  .$$
  \end{enumerate}
  \end{lemma}

\begin{definition}
  If $n = q^t$, where $q$~is prime, and $0 \le i < t$, we let
$\phi_i \colon q^i (\Z_n)^p \to (\Z_q)^p$ be the
natural homomorphism with kernel $q^{i+1}(\Z_n)^p$.
  \end{definition}

\begin{lemma} \label{reducetoprime}
  Let $n = q^t$ where $q$~is prime, and $p \neq q$, and let
$G$ be any transitive group of degree~$p$ that contains
$\Z_p$.
  \begin{enumerate}
  \item \label{reducetoprime-contains}
  If $C$ is any  $G$-invariant subgroup of~$(\Z_n)^p$, define
$C_i = \phi_i \bigl( C \cap q^i (\Z_n)^p \bigr)$ for $0 \le
i < t$. Then $C_0 \subset C_1 \subset \cdots \subset
C_{t-1}$ is an increasing chain of $G$-invariant subgroups
of $(\Z_q)^p$.
  \item \label{reducetoprime-exists}
  If $G \le \AGL(1,p)$, or $G = A_n$, or $G = S_n$, then the
converse holds:
  For any increasing chain $C_0 \subset C_1 \subset \cdots
\subset C_{t-1} \subset (\Z_q)^p$ of $G$-invariant subgroups
of~$(\Z_q)^p$, there is a subgroup of~$(\Z_n)^p$, such that
  $\phi_i \bigl( C \cap q^i (\Z_n)^p \bigr) = C_i$, for $0
\le i < t$.
  \item \label{reducetoprime-unique}
  Each $G$-invariant subgroup of~$(\Z_n)^p$ is uniquely
determined by the corresponding chain
  $C_0 \subset C_1 \subset \cdots \subset C_{t-1}$ of
$G$-invariant subgroups of $(\Z_q)^p$.
  \end{enumerate}
  \end{lemma}

\begin{proof} (\ref{reducetoprime-contains}) This follows
from the observation that, for any $c \in C \cap
q^i(\Z_n)^p$, we have $qc \in q^{i+1}(\Z_n)^p$ and
$\phi_i(c) = \phi_{i+1}(qc)$.

(\ref{reducetoprime-unique}) Suppose there is
a code~$C'$, such that $C'_i = C_i$ for each~$i$. Let $M = C
\cap q(\Z_n)^p$. By induction on~$t$, we may assume that $C'
\cap q(\Z_n)^p = M$. Consider the composite homomorphism:
  $$ C_0
  \cong \frac{C}{M}
  \hookrightarrow \frac{C'+q(\Z_n)^p}{M}
  \to \frac{C'+q(\Z_n)^p}{C'}
  \cong \frac{q(\Z_n)^p}{M} .$$
  If $C \neq C'$, then this homomorphism is nontrivial, so
$C_0$ and $q(\Z_n)^p/M$ have a composition factor in
common. Because
  $$ M \subset M + q^{t-1}(\Z_n)^p
  \subset M + q^{t-2}(\Z_n)^p
  \cdots
  \subset M + q(\Z_n)^p $$
  is an increasing chain of $G$-submodules with quotients
  $$ (\Z_q)^p/C_{t-1}, (\Z_q)^p/C_{t-2},
\ldots,(\Z_q)^p/C_1 ,$$
  we conclude that $C_0$ has a composition factor in common
with $(\Z_q)^p/C_i$, for some $i \ge 1$. This is
impossible, because $C_0 \subset C_i$, and the
representation of~$G$ on $(\Z_q)^p$ is multiplicity
free. (In fact, the restriction to the subgroup $\Z_p$ is
multiplicity free, because there are $p$~distinct $p$th
roots of unity in an appropriate extension of~$\Fq$.)

(\ref{reducetoprime-exists})
  It suffices to show, for each~$i$, that there is a
$G$-invariant subgroup~$\hat C_i$ of $(\Z_n)^p$, such that
$\hat C_j = C_i$ for $j \le i$ and $\hat C_j = 0$ for $j >
i$. (For then we simply let $C = \langle \hat C_0,
\ldots,\hat C_{t-1}\rangle$.) Thus, we may assume that, for
some~$i$, we have $C_0 = C_1 = \cdots = C_i$ and $C_{i+1} =
C_{i+2} = \cdots =  C_{t-1} = 0$.
  Furthermore, may assume that $i = t-1$ (because $C' =
q^{k}C$ satisfies $C'_j = 0$ for $j \ge t-k$).

  If $C_i$ is the repetition code, let
   $C$ be the repetition code in~$(\Z_n)^p$. If $C_i$ is the
dual of the repetition code, then let $C$ be the dual of the
repetition code in~$(\Z_n)^p$; that is,
  $$ C = \left\{\, (z_1,\ldots,z_p) \in (\Z_n)^p : \sum_{i=1}^p
z_i \equiv 0 \pmod{n} \right\}.$$
  Thus, we may now assume that $C_i$ is neither the
repetition code nor its dual. Then, from
Lemma~\ref{modp-allgroups} and the assumption on~$G$, we see
that $G \le \AGL(1,p)$.
  Therefore $\Z_p \triangleleft G$, so, from
uniqueness~(\ref{reducetoprime-unique}), we see that every
$\Z_p$-invariant subgroup of $(\Z_n)^p$ is $G$-invariant.
Thus, we may assume that $G = \Z_p$.

In this case, the desired conclusion is a special case of
\cite[Thm.~37.4, p.~156]{Burrow}, but we give an explicit
construction.
  Let $f(x) \in \Fq[x]$ be the monic generating polynomial
for~$C_i$. Because $f(x)$ is a divisor of $x^p - 1$, and
$x^p - 1$ has no repeated roots, we know, from Hensel's
Lemma \cite[36.5, p.~145]{Burrow}, that there is a monic
polynomial $g(x) \in \Z_n[x]$, such that $g(x) \equiv f(x)
\pmod{q}$, $\deg (g) = \deg(f)$, and $g$ is a divisor of
$x^p - 1$ in $\Z_n[x]$. Now let $C$ be the ideal of $\Z_n[x]
/ (x^p - 1)$ generated by $g(x)$.
  \end{proof}

\begin{remark} \label{goodp-K-notPSL}
  In applying Lemma~\ref{reducetoprime} to the study
of subgroups of~$S_{p^2}$, one is interested only in the
case where $n$~is not prime and there is a subgroup~$L$
of~$S_p$, such that $n$ is a divisor of $|N_{S_p}(L)/L|$.
Note that $3^2 \nmid 11-1$, $2^2 \nmid 23-1$, and neither
$11$ nor~$23$ can be written in the form $(q^d - 1)/(q-1)$
for a prime-power~$q$. Therefore, we see from
Lemma~\ref{modp-allgroups} that if $G = \PSL(2,11)$,
$M_{11}$, or~$M_{23}$, then, in the cases of interest, $C_i$
must be either the repetition code or its dual. Thus, the
proof of Lemma~\ref{reducetoprime} is valid in these cases. It
is only when $\PSL(d,q) \le G \le \PGammaL(d,q)$
that the possible choices of~$K$
in Proposition~\ref{wreath-prop}(\ref{wreath-prop-K}) have not yet
been completely classified.
  \end{remark}

\begin{lemma}[{\cite{BardoeSin, Klemm, Mortimer}}]
  \label{modp-allgroups}
  Let
  \begin{itemize}
  \item[$\bullet$] $p$~and~$r$ be prime;
  \item[$\bullet$] $G$ be a transitive subgroup
of~$S_p$ that contains $\Z_p$, and
  \item[$\bullet$] $C$ be a nontrivial cyclic code
over~$\Z_r$ that admits~$G$ as a group of permutation
automorphisms.
  \end{itemize}
  If $C$ is neither the repetition code nor its dual, then
either
  \begin{enumerate}
  \item $\Z_p \le G \le \AGL(1,p)$, and $C$ is described
in Lemma~\ref{subaffine-lem} below; or
   \item $G = \PSL(2,11)$, $p = 11$, $r = 3$, and $C$ is
either the $(11,6)$ ternary Golay code or its dual;
or
   \item $G = M_{23}$, $p = 23$, $r = 2$, and $C$ is either
the $(23,12)$ binary Golay code or its dual; or
  \item $\PSL(d,q) \le G \le \PGammaL(d,q)$, $p =
(q^d - 1)/(q-1)$, $q$~is a power of~$r$, and $C$ is
described in Theorem~\ref{PSLmodules} below.
  \end{enumerate}
  \end{lemma}

\begin{proof}
  From Theorem~\ref{degreep}, we know that there are only a
few possibilities for~$G$. In each case, the desired
conclusion is a known result.
  \begin{itemize}
  \item If $\Z_p \le G \le \AGL(1,p)$, see
Lemma~\ref{subaffine-lem}.
  \item If $G = A_p$ (or~$S_p$), see \cite[Beispiele
9(a)]{Klemm}.
  \item If $G = \PSL(2,11)$ and $p = 11$, see
\cite[(J)]{Mortimer}.
  \item If $G = M_{11}$ or~$M_{23}$ (and $p = 11$ or~$23$,
respectively), see \cite[Beispiele 9(bc)]{Klemm}.
  \item Suppose $\PSL(d,q) \le G
\le \PGammaL(d,q)$ and $p = (q^d - 1)/(q-1)$.
  If $r \nmid q$, see \cite[\S3(C)]{Mortimer};
  if $r \mid q$, see Theorem~\ref{PSLmodules}.
  \end{itemize}
  \end{proof}

\subsubsection{Cyclic codes invariant under a given subgroup
of $\AGL(1,p)$}
  \label{subaffine-sect}
  Lemma~\ref{subaffine-lem} characterizes the cyclic codes of
prime length~$p$ that admit a given subgroup of
$\AGL(1,p)$ as permutation automorphisms. This result must be
well known, but the authors have been unable to locate it in
the literature.

\begin{lemma} \label{subaffine-lem}
  Let $f(x) \in \Fq[x]$ be the generating polynomial of a
cyclic code~$C$ of prime length~$p$ over~$\Fq$, and let $A$
be a subgroup of $\Fp^*$.
  \begin{enumerate}
  \item If $p \nmid q$, then $C$ is $A$-invariant if and only
if
  $f(x)$ is a factor of $f(x^a)$,
  for every $a \in A$.
  \item If $p \mid q$, then $C$ is $A$-invariant.
  \end{enumerate}
  \end{lemma}

\begin{remark}
  Suppose $p \nmid q$. For a given subgroup~$A$ of $\Fp^*$,
one can construct all of the $A$-invariant cyclic codes of
length~$p$ by the following method.

Let $\poly \subset \Fq[x]$ be the set of all
monic factors of the polynomial $x^p - 1$, and let $\polyirr$
be the subset consisting of those polynomials that are
irreducible over $\Fq$. Then $A$ acts on both $\poly$
and~$\polyirr$ by
  $$ f^a(x) = \gcd\bigl( f(x^a), x^p - 1 \bigr) .$$
  From the lemma, we see that $f(x)$ is the generating
polynomial of an $A$-invariant code if and only if $f^a =
f$, for every $a \in A$.

  If $F$ is any $A$-invariant subset~$F$ of~$\polyirr$ (that
is, if $F$ is any union of orbits of~$A$), then
  $ \prod_{f \in F} f(x) $
  is the generating polynomial of an $A$-invariant code, and
conversely, every $A$-invariant generating polynomial can be
constructed in this way.

In particular, the number of $A$-invariant cyclic codes is
$2^d$, where $d$ is the number of $A$-orbits on~$\polyirr$.
However, it is probably easier to calculate~$d$ by using the
formula
  $ d = 1 + |\Fp^* : \langle A, q \rangle|$.
  \end{remark}

\subsubsection{Codes that admit $\PSL(d,q)$}
   Bardoe and Sin \cite[Thm.~A]{BardoeSin} recently gave an
explicit description of the codes that admit $\PGL(d,q)$ as
a group of permutation automorphisms. (They
\cite[Thm.~C]{BardoeSin} also considered monomial
automorphisms, but we do not need the more general result.)
For the case of interest to us, where $(q^d-1)/(q-1)$ is
prime, we know that $\gcd(q-1,d) = 1$, so the natural
embedding of $\PSL(d,q)$ into $\PGL(d,q)$ is an isomorphism.
Therefore, the codes described in~\cite{BardoeSin} are
precisely the codes that admit $\PSL(d,q)$ as a group of
permutations.

Furthermore, the results of Bardoe and Sin yield an explicit
description of the image of each code under the Frobenius
automorphism (cf.~\cite[Thm.~A(b)]{BardoeSin}), so the
results generalize easily to any subgroup~$G$ of
$\PGammaL(d,q)$ that contains $\PSL(d,q)$. After some
necessary definitions, we state this slightly more general
version of \cite[Thm.~A]{BardoeSin}.

\begin{definition}
  Suppose $r$ is a prime number, $q = r^t$, and $p = (q^d -
1)/(q-1)$ is prime. Let $c$ be a divisor of~$t$.

Let $\BSset{c}$ denote the set of $t$-tuples
$(s_0,s_1,\ldots,s_{t-1})$ of integers satisfying (for $j =
0,1,\ldots, t-1$, and with subscripts read modulo~$t$):
  \begin{enumerate}
  \item $1 \le s_j \le d-1$;
  \item $0 \le r s_{j+1} - s_j \le (r-1) d$; and
  \item $s_{j+c} = s_j$.
  \end{enumerate}
  Let $\BSset{c}$ be partially ordered in the natural way:
  $(s_0',\ldots,s_{t-1}') \le (s_0,\ldots,s_{t-1})$
  if and only if $s_j' \le s_j$ for all~$j$.

Let $\BSset{c}_0 = \BSset{c} \cup \{(0,0,\ldots,0)\}$, and
extend the partial order on~$\BSset{c}$ to~$\BSset{c}_0$, by
making $(0,0,\ldots,0)$ incomparable to all other elements.
  \end{definition}

\begin{definition}
  A monomial $X = \prod_{i=1}^{d} X_i^{b_i} \in
\Z_r[X_1,X_2,\ldots,X_d]$ is a \emph{basis monomial} if
  \begin{itemize}
  \item $0 \le b_i < q$, for $i = 1,\ldots,d$;
  \item $\deg(X) = \sum_{i=1}^{d} b_i$ is divisible by~$q-1$;
and
  \item  $X \neq X_{1}^{q-1} X_{2}^{q-1} \cdots X_{d}^{q-1}$.
  \end{itemize}
  \end{definition}

\begin{definition} ({\cite[\S3.2]{BardoeSin}})
  Let $X = \prod_{i=1}^{d} X_i^{b_i}$ be a basis monomial.
  For each $e \in \{0,1,\ldots,t-1\}$, let
  $$ \deg^e(X) = \sum_{i=1}^{d} \phi^e(b_i) ,$$
  where $\phi$ is the permutation on $\{0,1,\ldots,q-1\}$
defined by
  $\phi(k) = rk + (1-q) \lfloor rk/q \rfloor$.
  (In other words, if we write $k = \sum_{j=0}^{t-1} a_j r^j$
as a $t$-digit number in base~$r$, then $\phi(k) = a_{t-1} +
\sum_{j=1}^{t-1} a_{j-1} r^{j}$ is the $t$-digit number
obtained by rotating the $t$~digits of~$k$, including the
leading~$0$'s.)

Define
  $$s(X) = \frac{1}{q-1} \bigl( \deg^0(X), \deg^1(X), \ldots,
\deg^{t-1}(X) \bigr) .$$
  Then $s(X) \in \BSset{t}_0$.
  \end{definition}

\begin{definition}
  Any basis monomial~$X$ defines an $\Fq$-valued
function~$f_X$ on the vector space~$\Fq^{d}$. Because
$\deg(X)$ is divisible by~$q-1$, we have $f_X(v) =
f_x(\lambda v)$, for every $\lambda \in \Fq^*$ and $v \in
\Fq^{d}$, so $f_X$ factors through to a well-defined
function~$\bar f_X$ on the projective space~$\P^{d-1}(\Fq)$.
  \end{definition}

\begin{theorem}[{cf.\ Bardoe-Sin \cite{BardoeSin}}]
\label{PSLmodules} 
  Suppose $r$ is a prime number, $q = r^t$, $p = (q^d -
1)/(q-1)$ is prime, and $\PSL(d,q) \le G \le
\PGammaL(d,q)$. 

For any ideal~$\mathcal{I}$ of the partially ordered set
$\BSset{t}_0$, let
  $ M_{\mathcal{I}} \subset \Z_r[\P^{d-1}(\Fq)]$
  be the span over~$\Z_r$ of the functions~$\bar f_X$, for
all basis monomials~$X$, such that $s(X) \in
\mathcal{I}$.
  Then $M_{\mathcal{I}}$ is $G$-invariant.

Conversely, for each $G$-invariant subspace~$M$, there is a
unique ideal~$\mathcal{I}$ of~$\BSset{t}_0$, such that $M =
M_{\mathcal{I}}$.
  \end{theorem}
  
\noindent  \vrule  height 6pt width \textwidth \break
{\normalfont \emph{Note from the authors:} The Bardoe-Sin Theorem does not 
directly apply to the case where $G$ properly contains $\PSL(d,q)$,
and 
Theorem~12 in
the published version of our paper made an incorrect statement about this
general case.
For an accurate discussion of the codes invariant under a subgroup 
of $\mathop{\mathrm{P \Gamma L}}(d,q)$, see
Proposition~1.4, Theorem~1.5, and Section~3 of
[J.\,D.\,Dixon and A.\,E.\,Zalesski, 
Finite imprimitive linear groups of prime degree,
\emph{J.~Algebra} 276 (2004), no.~1, 340--370].
Corollary~3.8 of that paper tells us every $\mathop{\mathrm{P S L}}(d,q)$-invariant code is also $\mathop{\mathrm{P \Gamma L}}(d,q)$-invariant, which is not what our 
published statement of Theorem~12 indicated.

We thank Primo\v z Poto\v nik for both pointing out our error and explaining where to find the correct result in the literature.
}
\\ \vrule  height 6pt width \textwidth

\subsection{Crossed homomorphisms}
   \label{crossed-homo}

\begin{theorem}
  Let
  \begin{itemize}
  \item $p$ be a prime;
  \item $H$ be either $A_p$, $S_p$, or subgroup of~$\AGL(1,p)$
that contains~$\Z_p$;
  \item $n$ be a natural number, such that either $n = 2$ or
$n \mid p-1$ or $n \mid m$, where $m$ satisfies $p =
(r^{d^m}-1)/(r^d-1)$ for some prime~$r$ and natural
number~$d$;
  \item $K$ be an $H$-invariant subgroup of $(\Z_n)^p$; and
  \item $\phi \colon H \to (\Z_n)^p/K$ be a crossed
homomorphism.
  \end{itemize}
  Then $\phi$ is cohomologous to a homomorphism from~$H$ to
$C_0/(K \cap C_0)$, where  $C_0$ is the repetition code in
$(\Z_n)^p$.
  \end{theorem}

\begin{remark}
  The conclusion of the theorem can be stated more concretely:
  If $\phi$ is not cohomologous to~$0$, then either
  \begin{enumerate}
  \item $H \le \AGL(1,p)$, and there is some $c \in \Z_n$,
and some generator~$h$ of~$H/\Z_p$, such that
$|h|(c,c,\ldots,c) \in K$ and, after replacing~$\phi$ by a
cohomologous cocycle, we have $\phi(h^a ,z) =
a(c,c,\ldots,c)$, for $a \in \Z$ and $z \in \Z_p$; or
  \item $H = S_p$, $n$~is even, and there is some $c \in
\Z_n$, such that $(2c,2c,\ldots,2c) \in K$ and, after
replacing~$\phi$ by a cohomologous cocycle, we have
  $$ \phi(h) =
  \cases
   {
   \hfil 0 + K &\hbox{if $g \in A_p$} \cr
    (c,c,\ldots,c) + K  & \hbox{if $g \notin A_p$}
  }$$
  \end{enumerate}
  \end{remark}

\begin{proof}
  Let $V = (\Z_n)^p/K$, let, and let $C_0^\perp$ be its dual.

  Because $\gcd(p,n) = 1$, we know that every element of
$C_V(\Z_p)$ has a representative in $C_{(\Z_n)^p}(\Z_p)$
(cf.~\cite[Thm.~5.2.3, p.~177]{Gorenstein}). Therefore
  $$C_V(H) \subset C_V(\Z_p)
  = \frac{C_{(\Z_n)^p}(\Z_p) + K}{K}
  = \frac{C_0 + K}{K}
  \cong \frac{C_0}{K \cap C_0} .$$
  Thus, it suffices to show that, after replacing $\phi$
by a cohomologous crossed homomorphism, we have $\phi(H)
\subset C_V(\Z_p)$.

  \begin{Case}
  Assume $H \le\AGL(1,p)$.
  \end{Case}
  Because $\gcd(p,n) = 1$, we know that $H^1(\Z_p,V) = 0$
\cite[Cor.~12.2.7, p.~237]{Cart-Eil}. Therefore, replacing
$\phi$ by a cohomologous cocycle, we may assume that
$\phi(\Z_p) = 0$. Because $\phi$ is a crossed
homomorphism, this implies that $\phi(H) \subset C_V(\Z_p)$.

  \begin{Case} \label{Ap-coho}
  Assume $H = A_p$.
  \end{Case}
  Assume that $n$ is prime. From Lemma~\ref{modp-allgroups},
we know that $(\Z_n)^p/K$ is either $C_0$ or~$C_0^\perp$.
Because $A_p$ is perfect (or $p = 3$, in which case $A_p =
\Z_p$), we know that $H^1(A_p,C_0) = 0$. From
\cite[Lem.~1]{KleschevPremet}, we know that
$H^1(A_p,C_0^\perp) = 0$.

Let $m$ be a divisor of~$n$, such that $n/m$ is prime. By
induction on~$n$, we may assume that $\phi$ is cohomologous
to a crossed homomorphism into $mV$. Then the preceding
paragraph implies that $\phi$ is cohomologous to~$0$.

  \begin{Case}
  Assume $H = S_p$.
  \end{Case}
  From Case~\ref{Ap-coho}, we may assume, after replacing
$\phi$ by a cohomologous crossed homomorphism, that
$\phi(A_p) = 0$. Therefore, $\phi(S_p) \subset C_V(A_p) =
C_V(\Z_p)$.
  \end{proof}

\begin{remark}
  To complete the classification of transitive subgroups
of~$S_{p^2}$, the following problems remain:
  \begin{itemize}
  \item For $\PSL(d,q) \le G \le \PGammaL(d,q)$, extend the Bardoe-Sin
Theorem \ref{PSLmodules} from a classification of subgroups modulo a prime to
a classification modulo a prime-power.
  \item Calculate $H^1(H,V)$ for $H =
\PSL(2,11)$ (with $p = 11$), $M_{11}$, $M_{23}$ and
$\PSL(d,q)$.
  \item For each nontrivial cohomology class, find an
explicit crossed homomorphism to represent it.
  \end{itemize}
  \end{remark}

\section{Applications} \label{ApplicationSection}

\subsection{The Cayley Isomorphism problem}

Let $H$ be a set, and $E\subseteq 2^H\cup 2^{2^H}\cup\ldots$.  We say
that the ordered pair $X = (H,E)$ is a {\it combinatorial object}.  We call
$H$ the {\it vertex set} and $E$ the {\it edge set}.  If $E\subseteq 2^H$,
then $X$ is a {\it hypergraph}.  An {\it isomorphism} between two
combinatorial objects $X = (H,E)$ and $Y = (H',E')$ is a bijection
$\delta\colon H\to H'$ such that $\delta(E) = E'$.  An {\it automorphism} of
a combinatorial object $X$ is an isomorphism from $X$ to itself.  Let $G$ be
a group and $X = (G,E)$ a combinatorial object.  Define $g_L\colon G\to G$ by
$g_L(h) = gh$ and let $G_L = \la g_L\colon g\in G\ra$.  Then $X$ is a {\it
Cayley object of $G$} if and only if $G_L\le\Aut(X)$.  A Cayley object $X$ of
$G$ is a {\it CI-object of $G$} if and only if whenever $X'$ is a Cayley
object of $G$ isomorphic to $X$, then some $\alpha\in\Aut(G)$ is an
isomorphism from $X$ to $X'$.  Similarly, $G$ is a {\it CI-group with respect
to ${\cal K}$} if and only if every Cayley object in the class of
combinatorial objects ${\cal K}$ is a CI-object of $G$, and a {\it CI-group}
if $G$ is a CI-group with respect to every class ${\cal K}$ of combinatorial
objects. It is known \cite{Pal} that $G$ is a CI-group if and only if $\vert
G\vert = 4$ or $G\cong\Z_n$, with $(n,\varphi(n)) = 1$.  Hence neither
$\Z_{p^2}$ nor $\Z_p^2$ is a CI-group unless $p = 2$, although $\Z_p^2$ is a
CI-group with respect to graphs \cite{G}.  We begin with a characterization of
when two Cayley objects of a $p$-group $G$ can be isomorphic provided their
automorphism groups share a common Sylow $p$-subgroup.

\begin{lemma}\label{iso}
  Let $X$ and $Y$ be Cayley objects of a $p$-group $G$, and $P$ a Sylow
$p$-subgroup of both $\Aut(X)$ and $\Aut(Y)$.  Then $X$ and $Y$ are isomorphic
if and only if there exists $\delta\in N_{S_G}(P)$ such that $\delta(X) = Y$.
  \end{lemma}

\begin{proof} ($\Rightarrow$) Let $\omega\colon X\to Y$ be an isomorphism.
Then $\omega^{-1} P\omega \subset \Aut(X)$, and $\omega^{-1}P\omega\le P_1$, a
Sylow $p$-subgroup of $\Aut(X)$.  Hence there exists $\beta\in\Aut(X)$ such
that $\beta^{-1}P_1\beta = P$, so that $\beta^{-1}\omega^{-1}P\omega\beta\le
P$, which means $\omega \beta \in N_{S_G}(P)$. Furthermore, $\omega
\beta \colon X\to Y$ is an isomorphism.
  \end{proof}

\begin{corollary}\label{circiso}
Let $X$ and $Y$ be Cayley objects of $\Z_{p^2}$, such that $P_i$ is a Sylow
$p$-subgroup of both $\Aut(X)$ and $\Aut(Y)$, for some $2\le i\le p - 1$.  Let
$\beta\in\F_p^*$ such that $\vert\beta\vert = p - 1$.  Then $X$ and
$Y$ are isomorphic if and only of they are isomorphic by $\alpha =
\hat{\beta}^j\gamma_{i+1}^k$, for some $1\le j\le p - 1$ and $1\le k
\le p$.
\end{corollary}

\begin{proof} ($\Rightarrow$) From Lemmas \ref{iso} and~\ref{Np2}, we know
that $X$ and $Y$ are isomorphic if and only if they are isomorphic by some
$\delta\in N_{S_{p^2}}(P_i) = \la N_{S_{p^2}}(\la\tau\ra),\gamma_{i+1}\ra$.
As $P_i\tl N_{S_{p^2}}(P_i)$ and $\vert N_{S_{p^2}}(P_i)/P_i\vert = (p-1)p$,
there are $(p-1)p$ cosets of $P_i$ in $N_{S_{p^2}}(P_i)$. As
$\hat{\beta},\gamma_{i+1}\not\in P_i$, these $(p - 1)p$ cosets are
$P_i\beta^j\gamma_{i+1}^k$, $1\le j\le p - 1$ and $1\le k\le p$. Hence
$\delta$ may be written in the form $\delta = g\alpha$, with $g\in P_i$ and
$\alpha = \hat{\beta}^j\gamma_{i+1}^k$. Then $\alpha \in N_{S_{p^2}}(P_i)$
and, because $\delta(X) = Y$ and $g\in\Aut(Y)$, we have $\alpha(X) = Y$.
  \end{proof}

\begin{corollary}
  Let $X$ and $Y$ be Cayley objects of $\Z_p^2$ with $\Pi_1$ a Sylow
$p$-subgroup of $\Aut(X)$ and $\Pi_2$ a Sylow $p$-subgroup of $Y$.  Let
$\alpha_1\in\Aut(\Z_p^2)$ such that $\alpha_1 \Pi_1 \alpha_1^{-1} = P_i'$ and
$\alpha_2\in\Aut(\Z_p^2)$ such that $\alpha_2 \Pi_2 \alpha_2^{-1} = P_i'$,
$1\le i\le p - 1$.  Let $\beta\in\F_p^*$ such that $\vert\beta\vert = p -
1$.  Then $X$ and $Y$ are isomorphic if and only if they are isomorphic by
$\alpha_1\bar{\beta}^j\tilde{\beta}^k\gamma_{i+1}^\ell\alpha_2^{-1}$, $1\le
j,k\le p - 1$, $1\le \ell\le p$.
  \end{corollary}

\begin{proof} Note that $P_i$ is a Sylow $p$-subgroup of both
$\Aut(\alpha_1(X))$ and $\Aut(\alpha_2(Y))$.  It follows then by arguments
analogous to those in Corollary \ref{circiso} that $\alpha_1(X)$ and
$\alpha_2(Y)$ are isomorphic if and only if they are isomorphic by some
$\omega\in N_{S_{p^2}}(P_i')$. The result follows.
  \end{proof}

We remark that the case $P = P_2'$ was considered in \cite{B1}.

\subsection{Automorphism groups of Cayley graphs of $\Z_p^2$}

Using Theorem \ref{main3} we can calculate the full automorphism group of any
vertex-transitive graph of order $p^2$. We actually will prove this result in
slightly more generality, determining all $2$-closed groups $G$ that contain
a regular subgroup isomorphic to $\Z_p^2$ (as was done in the previously
cited paper). We remark that Klin and P\"oschel \cite{KPb} have already
calculated the full automorphism groups of circulant graphs of order $p^k$
(that is, of Cayley graphs of $\Z_{p^k}$).

\begin{theorem}\label{autozp2}
  Let $G$ be a $2$-closed subgroup of $S_{p^2}$ such that $G$ contains the left
regular representation of $\Z_p^2$.
  \begin{enumerate}
  \item \label{autozp2-doubly}
  If $G$ is doubly transitive, then $G = S_{p^2}$.
  \item \label{autozp2-primsolv}
  If $G$ is simply primitive and solvable, then $G\le\AGL(2,p)$.
  \item \label{autozp2-primnonsolv}
  If $G$ is simply primitive and nonsolvable, then $G\le AGL(2,p)$ or $G = S_2\wr S_p$ in
its product action.
  \item \label{autozp2-imprimsolv}
  If $G$ is imprimitive, solvable, and has elementary abelian Sylow
$p$-subgroup, then either $G < \AGL(1,p)\times\AGL(1,p)$ or $G = S_3 \times
S_3$ (and $p = 3$).
  \item \label{autozp2-imprimnonsolv}
  If $G$ is imprimitive, nonsolvable, and has elementary abelian Sylow
$p$-subgroup, then either $G = S_p\times S_p$ or $G = S_p\times A$, where $A
< \AGL(1,p)$.
  \item \label{autozp2-bigSylow}
  If $G$ is imprimitive with Sylow $p$-subgroup of order at least $p^3$,
then $G = G_1\wr G_2$, where $G_1$ and $G_2$ are $2$-closed permutation
groups of degree $p$.
  \end{enumerate}
  \end{theorem}

\begin{proof} (\ref{autozp2-doubly}) If $G$ is doubly transitive, then clearly
$G = S_{p^2}$.

(\ref{autozp2-imprimnonsolv}) If $G$ is imprimitive, nonsolvable, and 
has elementary abelian Sylow
$p$-subgroup, then by Theorem \ref{main3}, we have that $G =
\{(\sigma,\tau)\in H\times N_{S_p}(K):f(\sigma)\in \tau K\}$, where $K,H\le
S_p$ and $f\colon H\to N_{S_p}(K)/K$ is a group homomorphism.

Let $\tau_1\in H$ be a
$p$-cycle and $\tau_2\in K$ be a $p$-cycle.  Then $(\tau_1,1_{S_p})\in G$
and $(1_{S_p},\tau_2)\in G$.  Furthermore, $G$ admits complete block
systems ${\cal B}_1$ and ${\cal B}_2$ of $p$ blocks of cardinality $p$
formed by the orbits of $\la(\tau_1,1_{S_p}))\ra$ and
$\la(1_{S_p},\tau_2)\ra$, respectively (because $G\le S_p\times S_p$).

  If both $\Fix{G}{{\cal B}_1} = \{(\delta,1_{S_p}):\delta\in\Ker(f)\}$ and
$\Fix{G}{{\cal B}_2} = \{(1_{S_p}),\gamma):\gamma\in K\}$ are solvable, then
$\Fix{G}{{\cal B}_1} \le\AGL(1,p)$ and $\Fix{G}{{\cal B}_2} \le\AGL(1,p)$.
Then $K\le \AGL(1,p)$ and $N_{S_p}(K) = \AGL(1,p)$ is solvable.  Hence both
$\Ker(f)$ and $f(H)$ are solvable so that $H$ is solvable.  Thus $G$ is
solvable, a contradiction.

We now know that either $\Fix{G}{{\cal B}_1}$ or $\Fix{G}{{\cal B}_2}$ is
nonsolvable.  We will show that if $\Fix{G}{{\cal B}_2}$ is doubly transitive
(which includes the nonsolvable case), then $G = H\times S_p$. The case where
$\Fix{G}{{\cal B}_1}$ is nonsolvable is handled in a similar fashion.

If $\Fix{G}{{\cal B}_2}$ is nonsolvable, then by Theorem \ref{norsylp}
$\Fix{G}{{\cal B}_2}\vert_B$ is doubly transitive for every $B\in{\cal
B}_2$.  Hence $\Stab_{\Fix{G}{{\cal B}_2}}(i,j)\not = 1$ for every
$(i,j)\in\Z_p^2$.  Define an equivalence relation $\equiv$ on $\Z_p^2$ by
$(i,j)\equiv (k,\ell)$ if and only if $\Stab_{\Fix{G}{{\cal B}_2}}(i,j) =
\Stab_{\Fix{G}{{\cal B}_2}}(k,\ell)$.  As $G\le S_p\times S_p$, there are $p$
equivalence classes of $\equiv$ and each equivalence class of $\equiv$
contains exactly one element from each block of ${\cal B}_2$.  As
$\Fix{G}{{\cal B}_2}\vert_B$ is doubly transitive, $\Stab_{\Fix{G}{{\cal
B}_2}}(i,j)\vert_B$ has two orbits for every $B\in{\cal B}_2$.  One orbit
consists of $\{(k,\ell)\}$, where $(k,\ell)\equiv (i,j)$ and the other
consisting of the remaining elements of the block $B'$ of ${\cal B}_2$ that
contains $(k,\ell)$.  Let $\Gamma$ be an orbital digraph of $G$ with
$((i,j),(k,\ell))\in E(\Gamma)$.  If $i = k$, then $\Gamma = pK_p$ (the union
of $p$ disjoint copies of $K_p$) and so $\Aut(\Gamma) = S_p\wr S_p$.  If
$i\not = j$, then, as $\Gamma$ is an orbital digraph, either $(i,j)$ is only
adjacent to $(k,\ell)$ or $(i,j)$ is adjacent to every element of $B'$ except
$(k,\ell)$.  In either case, it is straightforward to verify that
$\{(1_{S_p},\gamma):\gamma\in S_p\}\le\Aut(\Gamma)$.  As $G$ is the
intersection of the automorphism group of all orbital digraphs of $G$, we
have $K = S_p$, $N_{S_p}(K) = S_p$ and $f = 1$.  Thus $G = H\times K = H\times
S_p$ as required.  Thus either $H < \AGL(1,p)$ or $H$ is doubly transitive (as
$\AGL(1,p)$ is doubly transitive).  Analogous arguments will then show that
if $H$ is doubly transitive, then $H = S_p$.  Thus
(\ref{autozp2-imprimnonsolv}) follows.

(\ref{autozp2-imprimsolv}) If $G$ is imprimitive, solvable, and has 
elementary abelian Sylow
$p$-subgroup, then we may define $H$, $K$, and~$f$ as in Theorem~\ref{main3}.
Both $H$ and $N_{S_p}(K)$ are solvable, so that $K$ is solvable and, by
Theorem \ref{norsylp}, we have $H,K\le\AGL(1,p)$.  As $N_{S_p}(\AGL(1,p)) =
\AGL(1,p)$, we have that $G\le\AGL(1,p)\times\AGL(1,p)$.  As $\AGL(1,p)$ is
itself doubly transitive, if $G = \AGL(1,p)\times\AGL(1,p)$ then
$\Fix{G}{\cal B}\vert_B$ is doubly transitive for every complete block system
${\cal B}$ of $G$ and every block $B\in{\cal B}$.  It then follows by
arguments above that $\Fix{G}{\cal B}\cong S_p$, a contradiction unless $p =
3$.  If $p = 3$, then $\AGL(1,p)\times\AGL(1,p) = S_3\times S_3$, a group
listed in (\ref{autozp2-imprimsolv}). Thus (\ref{autozp2-imprimsolv}) follows.

(\ref{autozp2-primsolv}, \ref{autozp2-primnonsolv}) If $G$ is simply
primitive, then by Theorem \ref{main3}, $G$ has an elementary abelian Sylow
$p$-subgroup and either $G\le\AGL(2,p)$ or $G$ contains an imprimitive
subgroup $H$ of index $2$. If $G \le \AGL(2,p)$, then the result follows, so we may assume $G$ contains an imprimitive subgroup $H$ of
index $2$. Note that $G$ is solvable if and only if $H$ is solvable.

If $H$ is solvable, then $H$ has an elementary abelian Sylow $p$-subgroup,
and so $G$ has an elementary abelian Sylow $p$-subgroup.  Furthermore, $G$ is
solvable.  Let $N$ be a minimal normal subgroup of $G$.  Then $N$ is an
elementary abelian $q$-group for some prime $q$.  As $G$ is primitive, $N$ is
transitive, $q = p$ and $\vert N\vert = p^2$.  Thus $G\le
N_{S_{\Z_p\times\Z_p}}(N)\le \AGL(2,p)$ and (\ref{autozp2-primsolv}) follows.

If $H$ is nonsolvable, then by (\ref{autozp2-imprimnonsolv}) proven above and
the fact that if $H\le G$, then $\cl(H)\le\cl(G)$, we have that either
$H = S_p\times S_p$ or $H = A\times S_p$, with $A < \AGL(1,p)$.  It then
follows by \cite[Theorem 4.6A]{DM} that $G = S_2\wr S_p$ with the product
action. Thus (\ref{autozp2-primnonsolv}) follows.

(\ref{autozp2-bigSylow}) If $G$ has a Sylow $p$-subgroup $\Pi$ of order at
least $p^3$, then $\Pi$ admits a complete block system ${\cal B}$ of $p$
blocks of cardinality $p$.  Then $\vert\Fix{\Pi}{\cal B}\vert\ge p^2$ so that
$\Stab_{\Fix{\Pi}{\cal B}}(0,0)\not = 1$.  As $\Fix{\Pi}{\cal B}$ is a
$p$-group, we have that if $\gamma\in\Stab_{\Fix{\Pi}{\cal B}}(0,0)$, then
$\gamma$ fixes every point of the block of ${\cal B}$ that contains
$(0,0)$.  Define an equivalence relation $\equiv'$ on $\Z_p^2$ by
$(i,j)\equiv'(k,\ell)$ if and only if $\Stab_{\Fix{\Pi}{\cal B}}(i,j) =
\Stab_{\Fix{\Pi}{\cal B}}(k,\ell)$.  It follows by comments above and the
fact that $\Stab_{\Fix{\Pi}{\cal B}}(i,j) = \Stab_\Pi(0,0)$, that the
cardinality of each equivalence class of $\equiv'$ is a multiple of $p$.  It
is straightforward to verify that the equivalence classes of $\equiv'$ are
blocks of $\Pi$ so that each equivalence class of $\equiv'$ has order $p$.
Thus the equivalence classes of $\equiv'$ form the complete block system
${\cal B}$.  For convenience, we assume without loss of generality that
${\cal B} = \bigl\{\{(i,j):j\in\Z_p\}:i\in\Z_p \bigr\}$.

Let $\Gamma$ be an orbital digraph of $G$, with $\Pi'$ a Sylow $p$-subgroup
of $\Aut(\Gamma)$ that contains $\Pi$.  Then $\Pi'$ admits ${\cal B}$ as a
complete block system as well.  If $\Gamma$ is disconnected, then $\Gamma =
p\Gamma_2$, where $p\Gamma_2$ is the disjoint union of $p$ copies of the
directed graph $\Gamma_2$ so that $\Aut(\Gamma) = S_p\wr\Aut(\Gamma_2)$.
If $\Gamma$ is connected, let $\bigl( (i,j),(k,\ell) \bigr)\in E(\Gamma)$ such
that $i\not = k$.  Then $(i,j)\not\equiv'(k,\ell)$ so that there exists
$\gamma\in\Pi$ such that $\gamma(i,j) = (i,j)$ but $\gamma(k,\ell)\not =
(k,\ell)$.  Then $\gamma$ permutes the $p$ elements of
$\{(k,m):m\in\Z_p\}$ as a $p$-cycle.  We conclude that
$\bigl((i,j),(k,m)\bigr)\in E(\Gamma)$ for every $m\in\Z_p$.  As
$\Fix{\Pi}{\cal B}$ is semiregular, we have that $\bigl((i,n),(k,m)\bigr)\in
E(\Gamma)$ for every $n,m\in\Z_p$.  Thus $\Gamma = \Gamma_1\wr\Gamma_2$ where
$\Gamma_1$ and $\Gamma_2$ are digraphs of order $p$.  It follows by
\cite[Theorem 1]{Sab} that $\Aut(\Gamma) = \Aut(\Gamma_1)\wr\Aut(\Gamma_2)$
(although the cited theorem is stated only for graphs, it works as well for
digraphs).  As $\cl(G)$ is the intersection of the automorphism groups of all
orbital digraphs of $G$, we conclude that $G = G_1\wr G_2$ for $2$-closed
groups $G_1,G_2$ of degree $p$.  Thus (6) holds.
  \end{proof}

\begin{theorem}\label{autozp21}
  Let $G$ be a $2$-closed subgroup of $S_{p^2}$ that contains the left regular
representation of $\Z_{p^2}$.  Then one of the following is true:
  \begin{enumerate}
  \item \label{autozp21-Sp2}
  $G = S_{p^2}$,
  \item \label{autozp21-N(Z)}
  $G\le N_{S_{p^2}}\bigl((\Z_{p^2})_L\bigr)$,
  \item \label{autozp21-wr}
  $G = G_1\wr G_2$, where $G_1$ and $G_2$ are $2$-closed groups of
degree $p$.
  \end{enumerate}
  \end{theorem}

\begin{proof}  If $G$ is doubly transitive, then $G = S_{p^2}$.  Otherwise,
as $\Z_{p^2}$ is a Burnside group \cite[Theorem 25.3]{W}, $G$ is
imprimitive.  By Theorem \ref{main3}, either a Sylow $p$-subgroup of $G$ is
normal in $G$, or a Sylow $p$-subgroup is isomorphic to $\Z_p\wr\Z_p$.  By
arguments in Theorem \ref{autozp2}, if a Sylow $p$-subgroup of $G$ has order
at least $p^3$, then $G = G_1\wr G_2$, where $G_1$ and $G_2$ are $2$-closed
groups of degree $p$.  The result then follows.
  \end{proof}

\begin{definition}
  A Cayley digraph $\Gamma$ of a group $G$ is {\it normal} if the left regular
representation of $G$ is normal in $\Aut(\Gamma)$.
  \end{definition}

In \cite[Problem 3]{MYX}, Ming-Yao Xu posed the problem of determining all
nonnormal Cayley graphs of order $p^2$.  We are now in a position to solve
this problem.

\begin{corollary} \label{nonnormal}
  A Cayley digraph~$\Gamma$ of a group of order $p^2$ is nonnormal if
and only if $\Gamma$ is isomorphic to one of the following graphs.
  \begin{enumerate}
  \item \label{nonnormal-complete}
  $\Gamma = K_{p^2}$, $p\ge 3$, or $p = 2$ and $G = \Z_4$,
  \item \label{nonnormal-wr}
  $\Gamma = \Gamma_1\wr\Gamma_2$, where $\Gamma_1$ and $\Gamma_2$ are
Cayley digraphs of the cyclic group of order $p$, $p\ge 3$,
  \item \label{nonnormal-prod}
  $\Gamma$ is a Cayley digraph of $\Z_p^2$ but not $\Z_{p^2}$, $p \ge
5$, with connection set $S = \{(i,0),(0,j):i,j\in\Z_p\}$ or the complement of
this graph,
  \item \label{nonnormal-weird}
  $\Gamma$ is a Cayley digraph of $\Z_p^2$ but not $\Z_{p^2}$, $p\ge 5$,
whose connection set $S$ satisfies the following properties, where $H =
\{(0,i):i\in\Z_p\}$,
  \begin{enumerate}
  \item $H\cap S = \emptyset$ or $H\cap S = H - \{(0,0)\}$,
  \item for every coset $(a,0) + H\not = H$ of $H$, $\bigl((a,0) + H\bigr)\cap
S = (a,b) + H,\emptyset,\{(a,0)\}$, or $\bigl((a,0) + H\bigr) - \{(a,0)\}$.
  \end{enumerate}
  \end{enumerate}
  \end{corollary}

\begin{proof} Let $G = \Aut(\Gamma)$. Then $\Gamma$ is normal if
(\ref{autozp2-primsolv}), (\ref{autozp2-imprimsolv}) of Theorem \ref{autozp2}
hold or (\ref{autozp21-N(Z)}) of Theorem \ref{autozp21} hold.

  If either (\ref{autozp2-doubly}) of
Theorem \ref{autozp2} or (\ref{autozp21-Sp2}) of Theorem \ref{autozp21}
holds, then $\Aut(\Gamma) = S_{p^2}$ and $\Gamma$ is not normal unless $p =
2$, in which case the left regular representation of $\Z_2^2$ is a normal
subgroup of $S_4$ but the left regular representation of $\Z_4$ is not a
normal subgroup of $S_4$ and (\ref{nonnormal-complete}) follows.

  If (\ref{autozp2-bigSylow}) of Theorem \ref{autozp2} or (\ref{autozp21-wr})
of \ref{autozp21} holds, then $\Gamma = \Gamma_1\wr\Gamma_2$ where $\Gamma_1$
and $\Gamma_2$ are Cayley digraphs of $\Z_p$.  It is then straightforward to
verify, as $\Aut(\Gamma) = \Aut(\Gamma_1)\wr\Aut(\Gamma_2)$ that left regular
representations of $\Z_{p^2}$ and $\Z_p^2$ are not normal in $\Aut(\Gamma)$
unless $p = 2$.  Whence (\ref{nonnormal-wr}) holds.  We conclude that the
remaining nonnormal Cayley digraphs must be Cayley digraphs of $\Z_p^2$ but
not $\Z_{p^2}$.

If (\ref{autozp2-primnonsolv}) of Theorem \ref{autozp2} holds, then
(\ref{nonnormal-prod}) follows by \cite[Theorem 2.12]{MYX}.

  Finally, if (\ref{autozp2-imprimnonsolv}) of Theorem \ref{autozp2} holds,
$\Gamma$ will be nonnormal provided that $p\ge 5$.  Further, $\Aut(\Gamma)$
admits a complete block system ${\cal B}$ of $p$ blocks of cardinality $p$,
which we may assume (by replacing $\Gamma$ with its image under an
appropriate automorphism of $\Z_p^2$) that ${\cal B}$ is formed by the orbits
of $H_L$ and that $\Fix{\Aut(\Gamma)}{\cal B}\vert_B = S_p$ for every
$B\in{\cal B}$.  As $\Fix{\Aut(\Gamma)}{\cal B}\vert_B = S_p$, we have that
$\Gamma[H] = K_p$ or $\bar{K_p}$, the complete graph on $p$ vertices or its
complement.  Whence $H\cap S = \emptyset$ or $H\cap S = H - \{(0,0)\}$.
Define an equivalence relation $\equiv$ on $\Z_p^2$ by $(i,j)\equiv (k,\ell)$
if and only if $\Stab_{\Fix{\Aut(\Gamma)}{\cal B}}(i,j) =
\Stab_{\Fix{\Aut(\Gamma)}{\cal B}}(k,\ell)$.  It is then straightforward to
verify that there are $p$ equivalence classes of $\equiv$ and that these $p$
equivalence classes of $\equiv$ form a complete block system ${\cal C}$ of
$\Aut(\Gamma)$.  Again, if necessary, we replace $\Gamma$ with its image
under an appropriate automorphism of $\Z_p^2$ and assume that ${\cal B}$ is
formed by the orbits of $H_L$ and ${\cal C} =
\bigl\{\{(i,j):i\in\Z_p\}:j\in\Z_p \bigr\}$.  Let $a\in\Z_p^*$.  Then $(0,0)$
is adjacent to either:
  no vertex of $(a,0) + H$;
  every vertex of $(a,0) + H$;
  only the vertex of $(a,0)$ of $(a,0) + H$; or
  every vertex of $(a,0) + H$ except $(a,0)$.
  Thus (\ref{nonnormal-weird}) follows.

The converse is straightforward.
  \end{proof}






\end{article}
\end{document}